\documentclass[11pt,letter]{article}

\usepackage{amsmath}
\usepackage{amssymb}

         \usepackage{epsf,epsfig}
            \usepackage{epic,eepic}
           \usepackage{rotating}

\usepackage{graphics}
\usepackage{color}

\def\qed{\hfill$\Box$\par\medskip\par\relax}

\newcommand{\8}{\infty}
\newcommand{\N}{{\mathbb N}}
\newcommand{\IP}{{\mathbb P}}
\newcommand{\E}{{\mathbb E}}
\newcommand{\R}{{\mathbb R}}

\newcommand{\Nn}{N_n}
\newcommand{\Din}{\Delta_{i,n}}
\newcommand{\Djn}{\Delta_{j,n}}
\newcommand{\Sn}{S_n}

\newcommand{\si}{\sigma}
\newcommand{\tn}{\tau_n}

\newcommand{\Tin}{T_{i,n}}
\newcommand{\Tnn}{T_{n,n}}
\newcommand{\Tun}{T_{1,n}}

\newcommand{\cTo}{{\mathcal T^{\circ}}}
\newcommand{\cToo}{{\mathcal T^{\boxtimes}}}
\newcommand{\cTt}{{\mathcal T^{\vartriangle}}}
\newcommand{\cW}{{\mathcal W}}
\newcommand{\rac}{{\varnothing}}
\newcommand{\card}{{\tt card}\;}
\newcommand{\Trans}{{\tt Trans}}

\newcommand{\ZGW}{Z^{\rm GW}}
\newcommand{\SGW}{{\tt Surv}^{\rm GW}}
\newcommand{\cTGW}{{\mathcal T^{\rm GW}}}
\newcommand{\siGW}{\sigma^{\rm GW}}

\newcommand{\rhops}{\rho^{\varepsilon}}
\newcommand{\ep}{\varepsilon}

\newtheorem{thm}{Theorem}[section]
\newtheorem{prop}[thm]{Proposition}

\newtheorem{lem}[thm]{Lemma}

\newcommand{\cvlaw}{\stackrel{\rm law}{\longrightarrow}}

\newcommand{\cF}{{\mathcal F}}
\newcommand{\Ceps}{C^{\varepsilon}}
\newcommand{\eps}{\varepsilon}
\newcommand{\un}{\frac{1}{n}}

{\markboth{}{\centerline{COMETS, DELARUE AND SCHOTT}}

\begin{document}

\thispagestyle{empty}

\title{Information Transmission 
under Random
  Emission 
  Constraints}

\author{Francis Comets${}^{(a),}$\footnote{comets@math.jussieu.fr; http://www.proba.jussieu.fr/$\sim$comets} \hspace*{.05cm},
Fran{\c c}ois Delarue${}^{(b),}$\footnote{delarue@unice.fr; http://math.unice.fr/$\sim$delarue} \hspace*{.05cm} and
Ren{\'e} Schott${}^{(c),}$\footnote{schott@loria.fr; http://www.loria.fr/$\sim$schott}
 \\ \\
(a) Laboratoire de Probabilit{\'e}s et Mod{\`e}les
Al{\'e}atoires,\\ Universit{\'e} Paris Diderot-Paris 7, Case 7012, 75205 Paris Cedex 13, France.
\\
\\
(b) Laboratoire J.-A. Dieudonn\'e, Universit\'e de Nice Sophia-Antipolis, 
\\
Parc Valrose, 06108 Nice Cedex 02, France.
\\ \\
(c) IECL and LORIA, Universit{\'e} 
de Lorraine, \\ 54506 Vandoeuvre-l{\`e}s-Nancy, France.}
\date{}
\maketitle
\begin{abstract} 
We model the transmission of a message on the complete graph with $n$ vertices and limited resources. The vertices of the graph represent servers that may broadcast the message at random. Each server has a random emission capital that decreases at each emission. 
Quantities of interest are the number of servers that receive the information before the capital of  all the informed servers is exhausted and the exhaustion time.
We establish limit theorems (law of large numbers, central limit theorem and large deviation principle), as $n\to \8$, for the proportion of informed vertices before exhaustion and for the total duration. The analysis relies on a  construction of the transmission procedure as a dynamical selection of successful nodes in a Galton-Watson
tree with respect to the success epochs of the coupon collector problem. 

\end{abstract}

\footnotesize
\noindent{\bf Short Title:} Information Transmission


\normalsize

\section{Introduction}

Transmission of information and dissemination of viruses in computer
networks gave rise to many practical as well as theoretical
investigations over the two last decades (see \cite{Bac, DZ, J, KTMN, KS, KLLM}).

In this paper, we model the transmission of a message on the complete graph with $n$ vertices and limited ressources.  Every vertex represents a server, which can be in one of three states: inactive (it did not receive the message yet), active (it has already received it, and is able to transmit it),  exhausted
(it has already received it, but cannot transmit it anymore because it has exhausted its own capital of emissions).  Each server $S_i$
has a random emission capital $K_i$. The message is initially received  from outside by one server, which is then turned from the inactive state to the active one (if it has a positive emission capital) or exhausted 
(if its emission capital is 0), though the $n-1$ other servers are inactive. At each integer time, one of the active servers (say $S_{i}$) casts the message, it looses one unit of its own emission capital $K_{i}$, and it 
selects the target at random among the $n$ servers. If the target is inactive, it discovers the information, it becomes itself active or exhausted according to its own emission capital. If not, 
this broadcast is unsuccessful and nothing else happens. When an
active server exhausts its emission capital, it enters the exhausted
state. The transmission ends at a finite time $\tau_n$, which is at
most equal to $1$ plus the sum of all initial capitals. 

From a practical point, the graph may be thought as a wireless
network, the vertices of which are battery powered sensors with a
limited energy capacity. We refer the reader to \cite{Bac, DZ, J, KTMN} for applications of graph theory to the
performance evaluation of information transmission in wireless networks.
We mention that the transmission process  can be also interpreted as the busy period of a queue
when the probability that a new customer enters the queue decays linearly as the number of past arrivals increases and finally vanishes after the $n$th arrival.

Here we describe the asymptotic behavior of the proportion of informed vertices at the end of the process when $n$ tends to the infinity. The 
mathematical analysis relies on a twofold structure: a subtree of the Galton-Watson tree, which models the vertices reached by the emission procedure, and the success epochs of 
the coupon collector problem, which model the successful transmissions. Phrased in a probabilistic way, we propose a coupling of the transmission model as a marginal tree of the Galton-Watson tree, obtained by pruning some of the nodes according to the coupon collector problem. Such a coupling provides a direct interpretation of the scenarios when the network ceases broadcasting 
at the very beginning of the process:
basically, these scenarios correspond to the extinction event in the Galton-Watson tree. On the survival event, we manage to specify the first-order behavior (in $n$) of the exhaustion time $\tau_{n}$ and of the proportion of informed nodes. Under suitable integrability conditions on the distribution of the capital of a given vertex, the fluctuations of both the exhaustion time and the proportion of informed nodes are also investigated: a central limit theorem is proved under a square-integrability condition and a somewhat involved large deviation principle is established under an exponential-integrability condition. In particular, when the distribution of the capital of a given vertex is of finite expectation, the probability that all the servers be reached before exhaustion (also referred to as the probability of full transmission) converges to $0$ as $n$ tends to the infinity;  as a consequence of the large deviation principle we prove here, it decays exponentially fast when the capital has a finite exponential moment. In some cases when the distribution of the capital has a heavy tail, we prove that the limit of the probability of full transmission is different from $0$ and, in particular, may coincide with the entire probability of survival of the Galton-Watson tree.
We refer to \cite{BaumBillingsley, BGL, DNW, FGT, Kan} and Chapter 8 of \cite{SF} for 
specific results concerning the coupon collector problem.


The papers  \cite{MMM} and
\cite{KLLM} study closely related transmission models. Machado et al. \cite{MMM} consider the case where $K_i=2$ and prove
partial transmission results. Obviously, our approach extends this result, as constant capitals are a particular  case of random ones. 
A specific interest of random capitals consists in allowing $K_{i}$ to be $0$ with a non-trivial probability: as we shall see below, 
a quick stop of the transmission process then occurs with a 
positive probability, as the extinction event of the Galton-Watson tree.
As in  \cite{MMM}, Kurtz et al. \cite{KLLM} investigate the
case where the $K_i$'s are constant, but possibly larger than $2$, time running continuously. In
their model, there is one particle at each vertex of the
graph at time 0; one of them is active, the others are inactive. The active
particle begins to move as a continuous-time, rate 1, random walk on
the graph; as soon as any active particle visits an inactive one,
the latter becomes active and starts an independent random walk. Each
active particle dies at the instant it reaches a total of L jumps
(consecutive or not) without activating any particle. Each active
particle starts with L lives and looses one life unit whenever it jumps on
a vertex which has already been visited by the process. 
For  another similar model with simultaneous jumps in discrete time, the number of informed servers has fluctuations 
of order $n^{3/4 }$ \cite{Zh12}.

We also emphasize that the dynamics of the present model is very similar, except for the asynchronisation,  to the frog model on the complete graph with finite lifetimes: the earliest reference is \cite{AlLeMaMa06},
addressing the question of final coverage. 
In the frog model on ${\mathbb Z}^d$, there is a phase transition between almost-sure extinction and 
survival with positive probability according to the underlying death rate \cite{AlMaPo02}, and similarly for a time-continuous model \cite{KeSi06}.
Shape theorems are proved in \cite{AlMaPo02b},  and also in \cite{RaSi04} and \cite{KS} in 
continuous time. Fluctuations are Gaussian in one dimension \cite{CoQuRa07, CoQuRa09}, but unknown when $d\geq 2$.

The paper is organized as follows. The basic model is presented in Section 2 together with the main results. In Section 3, we provide an alternative construction based on a pruning procedure of the Galton-Watson tree. Law of large numbers and related fluctuation limit theorems are investigated in Section 4, including the case of heavy tails. The large deviation principle is established in Section 5.


\section{The model and main results}
\subsection{A Markovian definition of the dynamics}
\label{subsec:markov}
From a modeling point of view, we assume that the servers emit between consecutive integers. The global state of the whole system before and after emissions is thus described at integer times. At any time $t \in {\mathbb N}$, $N_{n}(t)$ denotes the number of servers which have already received the message and the so-called `total emission capital' $S_{n}(t)$ the number of available attempts that can be used to deliver the message to a server which has not received it yet. At time $0$, only one server detains the information. Its own capital, that is $S_{n}(0)$ with our notation, is a random variable, the distribution of which is denoted by $\mu$. 

The dynamics of the pair process $(N_{n}(t),S_{n}(t))_{t \in {\mathbb N}}$ are then assumed to be Markovian. Conditionally on the states up to time $t$, the values of $N_{n}(t+1)$ and $S_{n}(t+1)$ are then given by:
\begin{equation}
\label{eq:ddef1}
S_n(t+1)=  S_n(t) + \;\left\{ 
\begin{array}{l}
-1 \\ K(t+1)-1 
\end{array}\right. \;{\rm with \ probability\ } 
\;\left\{ 
\begin{array}{l}
\Nn(t)/n \\  1-\Nn(t)/n
\end{array} \right.,
\end{equation}
where  $K(t+1)$ is a random variable, with $\mu$ as distribution and independent of the past up until time $t$; respectively, in the above cases,
\begin{equation}
\label{eq:ddef2}
N_n(t+1)=  \;\left\{ 
\begin{array}{l}
 N_n(t) \\  N_n(t)+1 
\end{array} \right. {\rm accordingly}.
\end{equation}
The Markov chain is absorbed at $S_n=0$. 
From a practical point of view, the interpretation is the following. During the emission that occurs between times $t$ and $t+1$, one server is chosen at random among the $n$ ones; it is referred to as the `target'. If the target is a server that has already received the information, then the number of informed servers remains the same and the total emission capital decreases by one. Such a scenario happens with the conditional probability $N_{n}(t)/n$. If the target is a non-informed server, then the number of informed servers increases by one and the total emission capital increases by the own emission capital of the target which has just been activated. This happens with conditional probability $1-N_{n}(t)/n$ and $K(t+1)$ then denotes the  initial emission capital of the target activated between time $t$ and time $t+1$. It is worth noting that the process $(N_{n}(t))_{ t \in {\mathbb N}}$ is a Markov process itself, known as the `coupon collector process' in the standard probability literature. It describes the collect, with replacement, of $n$ equally likely coupons. 

With such a modeling, the variables $K(t),t \in {\mathbb N},$ are i.i.d., with $\mu$ as common distribution, $K(0)$ being equated with $S_{n}(0)$ and the sequences $(K(t))_{t \in {\mathbb N}}$ and $(N_{n}(t))_{t \in {\mathbb N}}$ being independent (we emphasize that, when the emission between $t$ and $t+1$ is a failure, the variable $K(t+1)$ has no role in the description of the dynamics of the pair process $(N_{n}(t),S_{n}(t))_{ t \in {\mathbb N}}$). By independence of the two sequences, the sum of all the capitals revealed up until time $t$ has the same law as $R(N_{n}(t))$ where 
\begin{equation}
\label{eq:4:11:R}
R(t) = \sum_{s=1}^t K(s), \quad t \in {\mathbb N} \setminus \{0\}. 
\end{equation}

The transmission process lasts for a duration $\tau_n$ which is
the first time $t$ when  the emission capital is equal to 0,
\begin{equation}
\label{eq:tau=min}
 \tn = \min \{ t \in {\mathbb N}: \Sn(t)=0\}.
\end{equation}
A natural question consists in determining whether the information will reach all servers, or a proportion of them only. 
We then define the event of full transmission,
$$
\Trans_n = \{ N_n(\tn)=n\},
$$
which occurs when 
all the servers finally receive the information. Then, three regimes of interest can be distinguished, according to 
$${\mathbb P}(\Trans_n) = 0, \; {\rm or} \;1 ,  \; {\rm or} \; \in (0,1),$$
which naturally correspond to different tail behavior of $K$. 
In all these cases,  one is interested in the large-$n$ asymptotics of 
$\tn$ and $ \Nn(\tn)$.

\medskip

Then, with the notation $s\wedge t=\min\{s,t\}$, the sequence $(S_n(t \wedge \tau_{n}), N_n(t \wedge \tau_{n}))_{t \in {\mathbb N}}$ is a Markov chain on ${\mathbb N}\times \{0,1,\ldots,n\}$ with a non-decreasing second component, and
absorption on the vertical axis. 
The harmonic equations for absorption probabilities are rather intricate for a general $K$, a natural route being to 
approximate the process by a differential equation. Here, our analysis will rely on a specific construction of the dynamics obtained by considering the Galton-Watson tree of reproduction law $\mu$.

\subsection{Construction as a labeled Galton-Watson tree}

In Section \ref{sec:LGW}, we  construct the information transmission process with a Galton-Watson
tree with degree $K \sim \mu$ and an independent coupon collector process with $n$ coupons. We only give a quick account here. For each $n$, we prune the tree using the events of the $n$-coupon collector. 
Given a realization of the Galton-Watson tree, we visit 
successively all the nodes starting from the root, we keep [resp., erase] the current node if a new coupon is obtained at that time [resp., if we do not collect a new coupon; then, the whole
subtree below the current node is deleted]. This results into a subtree of the original one, with size $N_n(\tn) \leq n$. 

Here is a precise statement. We denote by $Z_k^{\rm GW}$  the cardinality of the $k$-th generation of the Galton-Watson tree and by $\SGW$ the survival event $  \SGW=\{ Z_k^{\rm GW}\geq 1, \forall k \}$. It is well known that ${\mathbb P}(\SGW)>0$ if and only if $\E K > 1$ or ${\mathbb P}(K=1)=1$.

\begin{prop} \label{prop:lgw}
Let
$(\Omega, {\mathcal A}, {\mathbb P})$ be   a probability space 
where are defined: (i) a  Galton-Watson tree with offspring distribution $\mu$, (ii) for each integer $n$, a  coupon collector process with $n$ coupons, independent of the tree.

Then, there exist sequences $(S_n(t), t \in {\mathbb N})$ and $(N_n(t), t  
\in {\mathbb N})$ defined on this probability space with $N_n(0)=1, S_n(0) \sim \mu$, and such that $((N_n(t),S_n(t)),t \in {\mathbb N})$ is a Markov
chain with transitions as in (\ref{eq:ddef1}),(\ref{eq:ddef2}). 

Letting
\begin{equation*}
Z_{\rm tot}^{\rm GW} = \sum_{k \geq 0} Z_k^{\rm GW},
\end{equation*}
the full transmission event  is `asymptotically included' in the survival event 
 in the sense that
$$
\Trans_n \subset \left\{ Z_{\rm tot}^{\rm GW}
\geq n \right\} \quad {\rm where} \quad \left\{ Z_{\rm tot}^{\rm GW}\geq n
\right\}
 \searrow  \SGW\; {\rm as} \; n \nearrow \8,
$$
and conversely, the event of termination at time $o(n)$ ($o(n)$ standing for the Landau notation) in the information process converges to  
extinction in the Galton-Watson process: denoting by $\Delta$ the symmetric difference,
\begin{equation}
\label{eq:dea}
\lim_{\epsilon \to 0^+} \lim_{n \to \8} {\mathbb P}\left( \{\tn \geq n \epsilon\} \Delta 
\SGW \right)=0.
\end{equation}
\end{prop}

The construction is simple and natural, but it seems to be new and it turns out to be 
a powerful tool to analyze the information process. In all our results below, we consider this particular coupling of the information process with the Galton-Watson tree and the coupon collector.
\subsection{Limit for the information coverage and duration}	

For $\E K > 1$, define $\theta \in (0,\8)$ as the unique root of the equation 
\begin{equation}
  \label{eq:theta}
\frac{1-e^{-\theta}}{\theta}= \frac{1}{\E K}.
\end{equation}
Extend this definition by setting
$\theta=0 {\rm \  if}\; \E K \leq 1, \theta=\8 {\rm \ if \ } \E K =\8.$ 
The function $\E K \mapsto \theta$ is an increasing bijection from $[0, \8]$ to $[0,\8]$.
Let also
\begin{equation}
  \label{eq:p(theta)}
p=1-e^{-\theta} \in [0,1],
\end{equation}
and note from  (\ref{eq:theta}),  that when $\E K \in (1, \8)$, $p=1-e^{-\theta} \in (0,1)$ is the unique solution  of
\begin{equation}
  \label{eq:p}
p\; \E K= - \ln (1-p), 
\end{equation}
whereas $p=1$ when $\theta=\8$. Let ${\bf 1}_A$ denote the indicator function of the event $A$.
\begin{thm} \label{th:LLN}
Let $\E K \in (0, \8]$.  

(i) As $n \to \8$, 
$$ \tau_n/n \longrightarrow \theta \;{\bf 1}_{\SGW}$$
in probability, with $\theta$ defined by (\ref{eq:theta}), though
$$ N_n(\tau_n)/n \longrightarrow p\; {\bf 1}_{\SGW}$$
in probability, with $p=1-e^{-\theta}$.

(ii) If $\E K \leq 1$ ($\IP(K=1) <1$ by assumption), then
$$
\lim_{n \to \8} \left[ \tau_n  + 1 \right] =\lim_{n \to \8}  N_n(\tau_n)  = 
Z_{\rm tot}^{\rm GW}
$$
in probability. 
\end{thm}
As a straightforward consequence, we observe that $\IP(\Trans_n) \rightarrow 0$ as $n \rightarrow \infty$ whenever 
$\E K$ is finite, as $p$ is strictly less than $1$ in this case. In Theorem \ref{thm:decay} we shall prove that the rate of decay is exponential in the case when $K$ has a finite exponential moment.

\subsection{Gaussian fluctuations in the case of a light tail}
\begin{thm} \label{th:flgauss}
Assume $\E K >1$ and $\E K^2 < \8$.
Let $\si_K^2$ denote the variance of $K$. As $n \to \8$,  conditionally on $\SGW$, we have the convergence in law:
$$
n^{-1/2} \big( \tn - n \theta \big) \cvlaw {\mathcal N}(0, \si_\tau^2),
$$
with 
%
%
$\si_\tau^2 = [(1-p)\E K-1]^{-2} [p\si_K^2 +  (\E K)^2 \si_N(\theta)^2]$
and $\si_N(s)^2= e^{-s}(1-e^{-s})-se^{-2s}$.
Similarly,  conditionally on $\SGW$, we have the convergence in law:
$$
n^{-1/2} \big( N_n(\tn) - n p \big) \cvlaw {\mathcal N}(0, \si_p^2),
$$
with $\si_p^2=
[(1-p)\E K-1]^{-2} [p\si_K^2 e^{-2\theta}+  \si_N(\theta)^2].
$ 
\end{thm}
Theorem \ref{th:flgauss} extends results in \cite{KLLM, MMM}. We do not study the random fluctuations
any further such as the corrections to the law of large numbers, but we prefer to explore the resulting  regimes of full transmission. From Theorem \ref{th:LLN}, when ${\mathbb E}K = \infty$, $N_n(\tn)/n \to 1$ in probability,
 conditionally on survival: This leaves open the asymptotics of the probability $\IP(\Trans_n)$. Some cases are investigated in the next subsection. 
\subsection{Probability of full  transmission in the case of a heavy tail}
When $K$ has a fat tail,  fluctuations will serve the transmission process. The next result shows that the full  transmission
event converges to the survival event if the tail is very heavy.
\begin{prop}
\label{th:FT+CT}
If there exist $c >0$ and $\alpha \in (0,1)$ such that 
$\liminf_{\ell \rightarrow \infty} 
[\ell^{\alpha} \IP(K \geq \ell)] \geq c $, then, \quad
$$\IP(\Trans_n) 
\longrightarrow \IP(\SGW) \quad {\rm as } \, n \rightarrow \infty.$$
\end{prop}

Hence, fluctuations of sums of $K$'s variables play a crucial role for the occurrence of full transmission. 
When $\mu$ has a heavy tail, one of the servers activated during the transmission process has a 
large enough capital $K$ allowing it to contact all the other servers.
We now turn to the critical case,
when $K$ belongs to the domain of attraction of a stable law of index 1. 

\begin{thm}
\label{th:critique} 
Assume there exists $c >0$ such that 
$$
\IP(K \geq \ell) \sim \frac{c}{\ell}\;,\qquad \ell \to \8. $$ Then,  as $n \rightarrow \infty$,
$$
\IP(\Trans_n) \longrightarrow \IP(\SGW) \times \left\{
\begin{array}{ll}
\E ( \exp(- e^{-{\mathcal S}})) \;,  & c=1,\\
1  \;,  & c>1,\\
0  \;,  & c<1,
\end{array}
\right.
$$
where, in the case $c=1$, ${\mathcal S}= \lim_{n \to \8} n^{-1}(R(n)-n\ln n)$ 
is a  totally asymmetric Cauchy variable, the centering parameter of which depends on the distribution of $K$
(see \eqref{eq:4:11:R} for the definition of $R$).
\end{thm}
The particular form of the limit relates to the celebrated result of Erd\"os and R\'enyi \cite{ErRe61} for the
coupon collector, that is the time to collect all coupons has  a Gumbel limit law.
\subsection{Large deviations}

Following Chapter 1 in \cite{dembo:zeitouni}, we recall that a random sequence $(Z_n, n \geq 1)$, with values in some Polish space $\mathcal Z$, obeys a large deviation principle (LDP) with rate function $\mathcal I$ and speed $n$ if \\
\phantom{} \hspace{5pt} $\bullet$  ${\mathcal I}: {\mathcal Z} \to [0,\8]$ is lower semi-continuous;\\
\phantom{} \hspace{5pt} $\bullet$ for all closed subset $F\subset  {\mathcal Z}$,
$\limsup_{n \to \8} n^{-1} \ln \IP(Z_n \in F) \leq -\inf \{{\mathcal I}(z); z \in F\}$;\\
\phantom{} \hspace{5pt} $\bullet$ for all open subset $O \subset  {\mathcal Z}$,
$\liminf_{n \to \8} n^{-1} \ln \IP(Z_n \in O) \geq -\inf \{{\mathcal I}(z); z \in O\}$.
\medskip

 Assume all through this section that $K$ has exponentially small tails: 
$\E \exp (a_0 K) < \8 $
 for some $a_0>0.$
Then, by Cramer's theorem (e.g., Subsection 2.2.1 in \cite{dembo:zeitouni}), $R(n)/n$ obeys a LDP with rate $I$,
\begin{equation}  \label{eq:ld0}
I(u)=\sup\{ au-\ln \E \exp (a K); a\in \R\},\qquad u \in [0,\8),
\end{equation}
and $I(u)=\8$ for $u <0$. The function $I$ is convex, 
lower semicontinuous and has 
compact level sets $\{u \in \R : I(u) \leq c\}$, for $c \in [0,\infty)$. The domain of the rate function $I$ is defined as the set $\textrm{Dom}(I) $ of reals $u$ with finite $I(u)$. Here, we have 
$\textrm{Dom}(I) = [k_*, k^*] \bigcap {\mathbb R}$, with
$$
k_*=\min\{k : {\mathbb P}(K=k)>0\}, \qquad k^*=\sup\{k : {\mathbb P}(K=k)>0\}.
$$
We also recall the large deviations principle for the coupon collector process with $n$ coupons from Boucheron et al. \cite{BGL}, Dupuis et al. \cite{DNW}:
for all $t>0$, we have:
\begin{equation} \label{eq:ld1}
N_n(nt)/n \quad {\rm obeys \ a \ LDP\ with \ speed \ } n {\rm \ and \ rate\ function\ } J_t. 
\end{equation}
From (2.7) and Section 4.1 in \cite{DNW}, the rate function $J_t$  is convex, it is finite if and only if $r \in (0, t \wedge 1]$, with 
$\lim_{r \searrow 0} J_t(r)=\8$, and it has a finite limit as $r \to t \wedge 1$. For
$r  \in (0, t \wedge 1)$, it is given by
\begin{equation}
\label{eq:Jt(r)}
J_t(r) =
(1-r) \ln (1-r) + (t-r) \ln \rho(r,t) +t e^{-t\rho(r,t)},
\end{equation}
where $\rho(r,t)$ denotes the unique solution in $(0,\8)$ of
\begin{equation} \label{eq:rhodef}
\frac{1-e^{-t\rho}}{\rho}=r. 
\end{equation}
With these ingredients we define 
the function $\cF: \R_+^2 \to [0,\8]$  by
\begin{equation}
\label{eq:rateF}
\cF(r,t)=
\left\{
\begin{array}{lll}
 rI(t/r) +J_t(r) &{\rm if}& r>0,\\
\8 &{\rm if}& r=0, \ t >0,\\
0 &{\rm if}& r=t=0.
\end{array}
\right.
\end{equation}

\begin{figure}[!h] 
\includegraphics[width=18cm]{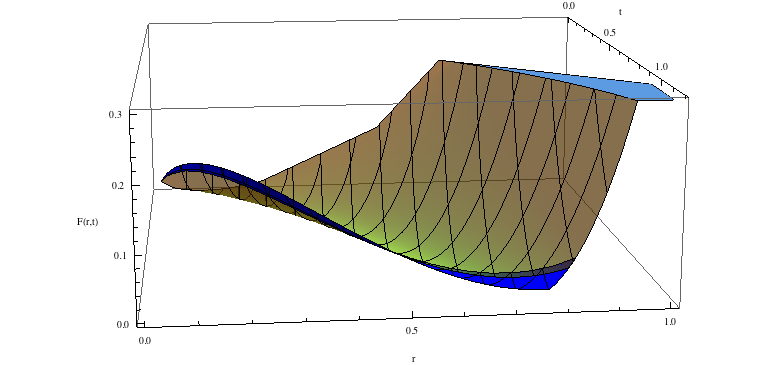} 
\caption{ Rate function $\cF$ for $K$ Poisson distributed with mean 1.4998. It vanishes at the origin and at  $(p, \theta) = (.5827,.8740)$, it is unbounded in neighborhoods of $(0,0)$ in its domain. 
For convenience, large values of $\cF$ are truncated, and the graph over the domain $r \leq .65 t$ is not shown. The dark blue strip corresponds to $.65 t \leq r \leq .67 t$, and the yellow part of the
graph to $ r \geq .67 t$.
 } 
\label{image_TAUX} 
\end{figure}

\begin{lem}
\label{lem:continuity}
(i) The function $\cF$,
is  lower semi-continuous on $\R_+^2$ with compact level sets $\{(r,t): \cF(r,t) \leq c\}$ for nonnegative $c$. Its domain is equal to
$${\rm Dom}({\mathcal F})=
\Big\{ (r,t) \in \R_+^2: 0 < r \leq t \wedge 1, (k_* \vee 1)r \leq t \leq k^* r\Big\} \cup \Big\{(0,0)\Big\}.
$$
It is continuous on 
${\rm Dom}({\mathcal F})\setminus \{(0,0)\}$. It is continuous at the origin if and only if $K$ is bounded.

(ii) Moreover, 
\begin{equation*}
\cF\bigl(1-e^{-t},t\bigr) = (1- e^{-t}) I \biggl( \frac{t}{1 - e^{-t}} 
\biggr),
\end{equation*}
and when $k_* \leq 1 < k^*$,  we have $\cF(1-e^{-t},t) \sim t I(1)$ as $t \rightarrow 0$. 

(iii) 
When $\E K >1$, the function $\cF$ is not convex, as it takes the value 0 at points $(p,\theta)$ and $(0,0)$, and is positive
elsewhere. 
When $\E K \leq 1$, $\cF$ is positive everywhere except at 0.
\end{lem}
The domain of $\cF$ might be rather degenerated. For example, in the 
Bernoulli case $K \in \{0,1\}$,
it reduces to the segment $\{(r,r), 0 \leq r \leq 1\}$. 
This function is an intricate combination of the rate functions of the coupon collector process and of the Galton-Watson process. This makes it an interesting rate function in its own. The shape of the graph of $\cF$ is shown in Figure \ref{image_TAUX} in the case $\E K >1$.

\begin{thm}
\label{th:ldp}
The sequence  $\big(n^{-1} (N_n(\tau_n), \tau_n); {n \ge 1 }\big)$ obeys a LDP with rate function $\cF$ and speed $n$.
\end{thm}
As a corollary, we obtain a variational formula for  the probability of full transmission.
\begin{thm}
\label{thm:decay}
The decay of the probability for all servers to be reached before exhaustion is exponential  and given by:
\begin{equation*}
\lim_{n \rightarrow +  \infty} n^{-1} \ln \IP ( \Trans_n ) 
= - 
\inf_{s \geq 0} \Bigl\{ I \bigl(\lambda(s) \bigr) + (\lambda(s) \!- \!1) \ln \bigl( 1\! - \!e^{-s}  \bigr) + \lambda(s) e^{-s}\Bigr\},
\end{equation*}
with $\lambda(s) = s/(1\!- \!e^{-s})$ for $s > 0,$ and $\lambda(0)=1$. The above right-hand side is negative.\end{thm}

\section{Construction from a labeled Galton-Watson tree} \label{sec:LGW}

%
%
%
Let $\cW= \cup_{k \geq 0} (\N^*)^k$ be the set of all finite words on the  alphabet $\{1,2,\ldots\}$. Its elements are of the form $w_1w_2\ldots w_k, w_i \in \N^*$ when $k \geq 1$, and,  for $k=0$, $(\N^*)^k$ reduces to the empty word  $\rac$, that we call the root. We then denote by $|w|=k$ the length of $w=w_{1}w_{2}\dots w_{k} \in \cW$ (with $\vert \rac \vert =0$.)  For $w, w' \in \cW$, we write $w < w'$
if: $|w| < |w'|$,  or $|w| = |w'|$ and $w {\leq}_{\rm lex} w'$ in the lexicographic order. We denote by $\preccurlyeq$
the usual predecessor relation in $\cW$, that is $w \preccurlyeq w'$ if $w$ is a prefix of $w'$.

Let $(K(w),w \in \cW)$ be a family of i.i.d. random variables on $\N= \{0,1,2,\ldots\}$ with common law $\mu$
(pay attention that the notation $K(t), t \in {\mathbb N}$, is also used in Subsection \ref{subsec:markov}; we
here use the same letter $K$, but a different index, no confusion being possible in the sequel). Assume $\mu(0)<1$ and $\mu(1)<1$ for a nontrivial setup.
The associated Galton-Watson tree $\cTGW$ is the set of $w \in \cW$ such that
$w = \rac$ or, for all $i=1,\ldots, |w|, w_i \leq K(v)$ 
with $v$ the predecessor of $w$ of length $i-1$ (in other words, given a parent $w'$ at the $(i-1)$th generation, that is 
$w'$ is a word of length $i-1$, the children of $w'$ are the words $w'1$, $\dots$, $w'K(w')$, of length $i$, obtained by concatenation). Denote by $\ZGW_k$ the size of the $k$th generation of this tree,
$\ZGW_k= \card  \{ w \in  \cTGW: |w|=k\}$, which is given by
$$
\ZGW_{k+1}= \sum_{v \in \cTGW, |v|=k} K(v), \quad \ZGW_0=1.
$$
Recalling that $\mu(1)<1$,
it is well known that the survival event $\SGW=\bigcap_k\{ \ZGW_k \geq 1\}= \{ \card  \cTGW = \8\}$
has complement probability 
$$
\siGW= 1- \IP(\SGW) =\left\{
\begin{array}{ll} 
=1 & {\rm if } \; \E K \leq 1,\\
<1 &  {\rm if } \; \E K > 1.
\end{array}
\right.
$$

On the same probability space, we consider an independent coupon collector process with $n$ images ($n\geq 1$):
Let  $\Din$, $i=1,\ldots n-1$, be independent,  geometrically distributed r.v.'s on
$\N^*$ with parameter $1-i/n$ (success probability) respectively. The success epochs are
$$
\Tun=0, \qquad \Tin= \sum_{j=1}^{i-1} \Djn, \qquad i=2,\ldots n,
$$
and the counting function is
$$\Nn(t)= \sum_{i=1}^{n} {\bf 1}_{\{\Tin \leq t\}}, \qquad t =0,1,\ldots
$$
In fact, $\Nn(t)$  represents the number of servers having received the information by time $t$ (note that $1 \leq N_n(t) \leq (t+1)\wedge n$).

For any fixed integer $n$,
with these two ingredients, we can define the transmission process together with the transmission time length $\tau_n$.
Let us start with an informal description. We browse a part of the Galton-Watson tree following the order $<$, and we paint the nodes in $\circ$ or in ${\triangle}$ according to the coupon collector process (success or failure); we only browse nodes which are in stand-by;
as soon as a node is painted in $\circ$, its number of children nodes in $\cTGW$ is revealed, and these children are put in stand-by. We then move to the next node in stand-by (next for $<$). 
The procedure runs until there are no nodes in stand-by anymore.

Here is a precise definition.
 Recursively for
$t=0,1,\ldots$, we construct $X(t) \in \cTGW$, and disjoint $\cTo(t), \cToo(t), \cTt(t) \subset \cTGW$  as follows
($X(t)$ encodes the vertex where  the $t$th tentative emission takes place, $\cTo(t)$ denotes the set of servers already informed by time $t$, $\cToo(t)$
is the set of tentative emissions scheduled but not yet performed at time $t$, $\cTt(t)$ is set of failed emissions, i.e. those
performed before time $t$ for which the target was already informed). 
Start with 
$$X(0)=\rac,\quad \cTo(0)=\{\rac\},\quad \cToo(0)=\{w \in \cTGW:  \vert w \vert =1, 1 \leq w_1 \leq K(\rac)\},\quad \cTt(0)=\emptyset.
$$
Here, and below, $\emptyset$ denotes the empty set and will not be confused with the root $\rac$ of the tree.
With the process $(X(t),\cTo(t), \cToo(t), \cTt(t))$ at time $t$, its value at the next step $t+1$ is defined by:
\begin{itemize}
\item If $\cToo(t)$ is nonempty, we let $X(t+1)$ be its first element, 
$$
X(t+1)= \inf\{ w \in \cToo(t)\}, \; {\rm denoted \ by} \; v
$$
to ease the notations, and we perform a test:
\begin{itemize}
\item If $N_n(t+1)=N_n(t)+1$, we define
\begin{equation}	\label{eq:deftrans1}
\begin{array}{l}
\cTo(t+1)=\cTo(t) \cup \{v\}, \\
 \cToo(t+1)=\big(\cToo(t) \setminus \{v\}\big) \cup \big\{v1, v2, \ldots, vK(v)\big\}, \\
\cTt(t+1)=\cTt(t),
\end{array}
\end{equation}
the notation $vk$ denoting the word of length $|v|+1$ obtained by concatenation.
\item 
 If $N_n(t+1)=N_n(t)$, we define
\begin{equation}	\label{eq:deftrans2}
\begin{array}{l}
\cTo(t+1)=\cTo(t), \\
 \cToo(t+1)=\cToo(t) \setminus \{v\}, \\
\cTt(t+1)=\cTt(t) \cup \{v\}.
\end{array}
\end{equation}
\end{itemize}
\item If $\cToo(t)$ is empty, we set $\tau_n=t$, and the construction is stopped (as well as the transmission). 
The set $\cTo(t)=\cTo(\tau_n)=\cTo(\8)$ is the set of servers finally informed. Note that $\tau_n \leq \Tnn$ is a.s. finite.
\end{itemize}
We observe that for all $t$, $\cTo(t)$ is a tree, as well as $\cTo(t) \cup \cToo(t)$. Moreover,  $\cToo(t) \cup \cTt(t)$
is a cutset of $\cTGW$ for its own graph structure.
\vspace{-40pt}

	\begin{center}
	\unitlength 1cm
	\begin{picture}(11,11)(-0.5,-0.5)%
	\linethickness{20cm}%

	\Thicklines
	\put(5,7.9894){\begin{rotate}{90}\ellipse{0.4266}{0.4266}\end{rotate}}
         \put(4.4,8.2){$\rac$}
	\color{black}
	\put(6.5079,6.2434){\begin{rotate}{90}\ellipse{0.4769}{0.4769}\end{rotate}}
%
%
	\put(8.7302,6.2169){\begin{rotate}{90}\ellipse{0.4791}{0.4791}\end{rotate}}
	\put(1.2434,6.2434){\begin{rotate}{90}\ellipse{0.4769}{0.4769}\end{rotate}}
	\put(8.7698,4.246){\begin{rotate}{90}\ellipse{0.5033}{0.5033}\end{rotate}}
	\put(9.7487,4.246){\begin{rotate}{90}\ellipse{0.5026}{0.5026}\end{rotate}}
	\put(7.7513,4.246){\begin{rotate}{90}\ellipse{0.535}{0.535}\end{rotate}}
	\put(2.2487,4.2328){\begin{rotate}{90}\ellipse{0.5026}{0.5026}\end{rotate}}
	\put(2.2222,2.2354){\begin{rotate}{90}\ellipse{0.5033}{0.5033}\end{rotate}}
	\path(9.5106,2.2487)(9.5106,1.7328)(8.9947,1.7328)(8.9947,2.2487)(9.5106,2.2487)
	\path(8.5053,2.2354)(8.5053,1.7328)(8.0026,1.7328)(8.0026,2.2354)(8.5053,2.2354)
	\path(7.2619,2.2354)(7.2619,1.746)(6.7725,1.746)(6.7725,2.2354)(7.2619,2.2354)
	\path(2.9762,0.7407)(2.9762,0.2381)(2.4735,0.2381)(2.4735,0.7407)(2.9762,0.7407)
	\path(1.9974,0.754)(1.9974,0.2381)(1.4815,0.2381)(1.4815,0.754)(1.9974,0.754)
	\path(3.5053,6.4947)(3.7566,6.0053)(3.2275,6.0053)(3.5053,6.4947)
	\path(.5026,4.4841)(.7407,3.9947)(.2381,3.9947)(.5026,4.4841)
	\path(9.5106,2.2487)(8.9947,1.7328)
	\path(9.5106,1.7328)(8.9947,2.2487)
	\path(8.5053,2.2354)(8.0026,1.7328)
	\path(7.2619,2.2354)(6.7725,1.746)
	\path(8.5053,1.7328)(8.0026,2.2354)
	\path(7.2619,1.746)(6.7725,2.2354)
	\path(2.9762,0.2381)(2.4735,0.7407)
	\path(1.9974,0.2381)(1.4815,0.754)
	\path(2.9762,0.7407)(2.4735,0.2381)
	\path(1.4815,0.2381)(1.9974,0.754)
	\path(8.5317,6.3757)(5.2052,7.9233)
	\path(6.3757,6.455)(5.1323,7.8307)
	\path(4.8545,7.8307)(3.5847,6.3492)
	\path(4.7751,7.9233)(1.4286,6.3889)
	\path(8.7566,5.9656)(8.7698,4.5106)
	\path(8.8889,6.0317)(9.6958,4.4841)
	\path(8.5582,6.0317)(7.8175,4.4974)
	\path(8.9153,4.0476)(9.2583,2.2487)
	\path(8.6376,4.0212)(8.2275,2.2322)
	\path(7.6323,3.9947)(7.01,2.2487)
	\path(1.336,6.0317)(2.1296,4.4577)
	\path(2.2619,3.9647)(2.2354,2.5)
	\path(2.3545,2.027)(2.7381,0.7407)
	\path(2.0635,2.0106)(1.7328,0.7672)
	\path(1.0979,6.0582)(0.5423,4.418)
	\path(3.4921,6.0053)(3.4921,4.2593)
	\path(3.4921,4.2593)(3.9947,2.25)
	\path(3.4921,4.2593)(3.0423,2.25)
	\path(0.6085,3.9683)(0.9921,2.3016)
	\path(0.3876,3.9815)(0.1058,2.25)
	\path(0.9921,2.25)(0.9921,0.7672)
	\thinlines
	\path(6.5868,1.57622)(6,1)
	\path(6.4568,1.7222)(6.7436,1.73)(6.7232,1.45)(6.4568,1.7222)
	\allinethickness{0.4pt}%
	\put(5.2,.8){$X(t+1)$}
	\put(3.4921,4.2593){\begin{rotate}{90}\ellipse*{0.15}{0.15}\end{rotate}}
	\put(0.9921,2.25){\begin{rotate}{90}\ellipse*{0.15}{0.15}\end{rotate}}
	\put(.1058,2.25){\begin{rotate}{90}\ellipse*{0.15}{0.15}\end{rotate}}
	\put(.9788,0.7407){\begin{rotate}{90}\ellipse*{0.15}{0.15}\end{rotate}}
	\put(3.0423,2.25){\begin{rotate}{90}\ellipse*{0.15}{0.15}\end{rotate}}
	\put(3.9947,2.25){\begin{rotate}{90}\ellipse*{0.15}{0.15}\end{rotate}}
%
	\put(-1.5,-.5){
	Figure 2. The Galton-Watson tree is represented up to the 4th generation.} 
	\put(-.7,-1.1){
	$\cTo(t)$, $\cToo(t)$ and $\cTt(t)$ and $X(t+1)$ are represented at time $t=10$.}
	\end{picture}%
	\end{center}

\vspace{30pt}

Figure 2 provides an example of construction of the sets $\cTo(t)$, $\cToo(t)$ and $\cTt(t)$, $t \in \{0,1,\dots,10\}$, according to the rules prescribed in
 \eqref{eq:deftrans1} and \eqref{eq:deftrans2}. At time $0$, the capital of the initial server is $K(\rac)=4$, so that the root of the tree has four children. Between times $0$ and $1$, the first emission is a success as the server which is revealed in the tree is painted in $\circ$. This server reads as the first child (starting from the left) at the first generation of the tree; it has two children, that is $K(\textrm{`1'})=2$, where `1' is here understood as a one-letter word. At time $1$, the total emission capital is thus $S(1)=4+2-1=5$, and $\cToo(2)$ contains 
 the three last children at the first generation and the two first children at the second generation. Since the second child at the first generation is painted in $\triangle$, the emission between times $1$ and $2$ fails, which means that the children of this node are not considered for the sequel of the construction. At time $2$, $S(2)=4+2-2=4$. Then, the emission between times $2$ and $3$ is a success but the server which is informed has no children (third child in the first generation) and the emission between times $3$ and $4$ is a success as well, with $K(\textrm{`4'})=3$. At time $4$, $\cToo(4)$ contains 5 nodes, all of them at the second generation of the tree. The node $X(5)$ is then the first child at the second generation: as it is painted in $\triangle$, the emission between times $4$ and $5$ fails.  
 And so on up until time $10$. Then, $\cToo(10)$ contains 5 nodes: two of them at the fourth generation and three of them at the third generation. With the lexicographic order, $X(11)= \textrm{`411'}$.

We now relate the above construction to the dynamical model for transmission.
Let $\card A$ denote the cardinality of a set $A$.
Consider a new, independent, i.i.d. sequence $(\bar K_i)_{i \geq 1}$ with law $\mu$, and define, for $i=1,\ldots n$, 
\begin{equation}
\label{def:ki}
K_i =
\left\{
\begin{array}{ll}
K(X(\Tin)) & {\rm if} \; i \leq \card  \cTo(\infty),\\
\bar K_i & {\rm if} \; i > \card  \cTo(\infty),
\end{array}
\right.
\end{equation}
and also
\begin{equation} 
\label{eq:S}
\Sn(t):= \sum_{i=1}^{\Nn(t)} K_i -t, \quad t \in {\mathbb N} \setminus \{0\}.
\end{equation}
Below, we also write
\begin{equation}
\label{eq:R:2}
\Sn(t) = R \left( \Nn(t) \right) -t, \quad t \in {\mathbb N}, \quad \textrm{with} \quad R(m) := \sum_{i=1}^m K_{i}, \quad m \in {\mathbb N} \setminus \{0\}. 
\end{equation}
Pay attention that the letter $R$ is also used in  \eqref{eq:4:11:R}, but as proved right below the two $R$'s have the same distribution. In the sequel, we will always refer to \eqref{eq:R:2} for the precise definition of $R$. 

By construction, we have
\begin{equation}\label{eq:identite}
N_n(t)=\card  \cTo(t), \quad S_n(t)=\card  \cToo(t).
\end{equation}
\begin{prop} \label{prop:eqmod}
The variables $(K_i)_{1 \leq i \leq n}$ are independent, identically distributed with law $\mu$, and 
$(K_i)_{1 \leq i \leq n}$ is independent of $(\Tin)_{1 \leq i \leq n}$. Moreover,
\begin{equation}
\label{eq:tau=min2}
 \tn = \min \{ t \in {\mathbb N}: \Sn(t)=0\}.
\end{equation}
\end{prop}
 {\bf Proof:} The formula for $\tn$ directly follows from \eqref{eq:identite}, and the fact that $\tau_n$ is finite.
We now investigate the distribution of the sequence $(K_{i},i \leq n)$.  
Below, we denote by ${\mathcal F}_{w} = \sigma(K(w'),w' \leq w)$ for $w \in {\mathcal W}$. 
On the event $A=\{(T_{1,n},\dots,T_{i,n})=(k_{1},\dots,k_{i}),X(k_{i}-1) = w,\tn\geq k_{i}\}$, 
$0=k_{1} < \dots < k_{i}$ and $w \in {\mathcal W}$, $\vert w \vert \leq k_{i}-1$, 
$X(T_{i,n})$ coincides with an ${\mathcal F}_{w}$-measurable r.v., denoted by 
$\chi$, which satisfies $w<\chi$ almost-surely
(this follows from the monotonicity of the browsing procedure). Similarly, all the variables $K(X(T_{j,n}))$, $1 \leq j \leq i-1$ coincide with ${\mathcal F}_{w}$-measurable r.v.'s on $A$.
Clearly, $K(\chi)$ is independent of 
${\mathcal F}_{w} \vee \sigma(T_{1,n},\dots,T_{n,n})$ and has $\mu$ as distribution, since the r.v.'s $(K(w'), w' \in {\mathcal W})$ are i.i.d and are independent of the success epochs $(T_{1,n},\dots,T_{n,n})$. Obviously, the event 
$A$ belongs to ${\mathcal F}_{w} \vee \sigma(T_{1,n},\dots,T_{n,n})$.
This proves that, for any bounded and measurable Borel function $\phi$, 
\begin{equation*}
{\mathbb E} \bigl[  \phi(K_{i}) {\mathbf 1}_{\{T_{i,n} \leq \tau_{n}\}} \vert (T_{j,n},1 \leq j \leq n),(K_{j},1\leq j \leq i-1) \bigr]
= {\mathbf 1}_{\{T_{i,n} \leq \tau_{n}\}} \int_{\mathbb N} \phi d \mu. 
\end{equation*}
On the event 
$\{(T_{1,n},\dots,T_{i,n})=(k_{1},\dots,k_{i}),\tau_{n} <k_{i}\}$, $K_{i}$ coincides with 
$\bar{K}_{i}$, which is obviously independent of $\sigma((T_{j,n},1 \leq j \leq n),(K_{j},1\leq j \leq i-1))$, so  that the above equality also holds with $T_{i,n} \leq \tau_{n}$ replaced by $T_{i,n} > \tau_{n}$.  \qed
Then the process we have constructed here corresponds to the description of the information
transmission process given in Subsection \ref{subsec:markov}. \medskip

\noindent
 {\bf Proof}
of Proposition \ref{prop:lgw}:
(\ref{eq:dea}) follows from Lemma \ref{lem:tau} below and $\{\tn \geq n\epsilon \} \subset 
\{Z_{\rm tot}^{\rm GW} \geq n\epsilon\}$.
The other  claims directly follow from the construction and Proposition \ref{prop:eqmod}. 
 \qed


\section{Proofs of law of large numbers and fluctuations}

In all the proofs, we use the following convention: for any discrete process $(V_{t},t \in {\mathbb N})$, $(V_{t},t \geq 0)$ stands for $(V_{\lfloor t \rfloor}, t \geq 0)$, where $\lfloor \cdot \rfloor$ is the floor function. We will also use the ceiling function, denoted by 
$\lceil \cdot \rceil$. For an interval $I$ in $\R$, define the Skorokhod space $D(I)$ as the space of 
c\`adl\`ag (right continuous left limited) functions from $I$ to $\R$. 

\subsection{Proofs of the law of large numbers}

\begin{lem}  \label{lem:PI}
We have the following convergence in law of sequences of processes on the Skorohod space (endowed with the standard J1 topology, keeping in mind that convergence for the J1 topology implies uniform convergence on compacts when the limit function is continuous): 

(i) On $D([0,1))$,
$$n^{-1/2} \big( T_{\lfloor nq \rfloor,n}- n \ln \frac{1}{1-q}\big)_{0 \leq q < 1} \stackrel{\rm law}{\longrightarrow}
\bigl[ B\big( \si_T(q)^2\big) \bigr]_{0 \leq q < 1} \quad \textrm{as} \ n \to \infty,$$
with $B$ a standard Brownian motion, and  
$$
 \si_T(q)^2= \frac{q}{1-q}+ \ln (1-q) >0.
$$ 

(ii) On $D(\R^+)$, 
$$
\bigl( n^{-1/2} \big( N_n(ns) -n(1-e^{-s}) \big) \bigr)_{s \geq 0} \cvlaw 
\bigl[
 B\big( \si_N(s)^2\big) \bigr]_{s \geq 0} \quad \textrm{as} \ n \to \infty,$$
with $B$ a standard Brownian motion, and  
$$
 \si_N(s)^2= e^{-s}(1-e^{-s})-se^{-2s} >0.
$$ 
Both limits are independent increments Gaussian processes with zero mean, and they are martingales. 
\end{lem}
 {\bf Proof:} Assertion (i) is a direct application of  the invariance principle for  triangular array of independent, but not i.d., square-integrable r.v.'s, see Dacunha-Castelle and Duflo \cite[Th\'eor\`eme 7.4.28]{dac:duflo}
or  Jacod and Shiryaev \cite[Chapter VII, Theorem 5.4]{jac:shy}. The variance is computed as the limit of a Riemann sum,
$$
 \si_T(q)^2= \lim_{n \to \8} \frac{1}{n} \sum_{i=1}^{qn} {\rm Var}(\Din)= \lim_{n \to \8} \frac{1}{n} \sum_{i=1}^{qn} \frac{i/n}{(1-i/n)^2}= \int_0^q \frac{y}{(1-y)^2}dy.
$$
Assertion (ii) follows from (i), using that $N_n(\cdot)$ and $T_{\lfloor n\cdot \rfloor,n}$ are reciprocal in a generalized sense. With 
$f(q)=-\ln(1-q),  f^{-1}(s)=1-e^{-s}$, we have $ \si_N(s)^2= \si_T(f^{-1}(s))^2 \times [f' \circ f^{-1}(s)]^{-2}$, see Billingsley \cite[Theorem 17.3]{billinglsey}. 
\qed
The next lemma is one of the key argument of the whole analysis. It shows that when the Galton-Watson tree is infinite transmission takes place on a macroscopic time level.
\begin{lem} \label{lem:tau}
There exists $\ep_0 >0$ such that for all $\ep \in (0,\ep_0)$, 
$$\lim_{n \to 
\8} \IP( \tau_n \geq n \ep, \SGW) =\IP(  \SGW) =1-\siGW.$$
\end{lem}
 {\bf Proof:} The claim being trivial for $\siGW=1$, we just need to consider the case when $\E K >1$.
Letting here $k=\lfloor \ln^2 n \rfloor$,
we estimate
\begin{equation}
  \label{eq:tau>}
\IP(\tau_n \leq n \ep, \SGW) \leq
  \IP(N_n(k) \leq  k) + 
  \IP(k \leq \tau_n \leq n \ep, \SGW) ,
\end{equation}
using that $\{ N_n(k)=k+1, \SGW\} \subset \{ \tau_n \geq k, \SGW\}$ which implies that 
$$\IP\left( N_n(k)\leq k, \SGW\right) \geq \IP \left( \tau_n < k, \SGW\right).$$
We start with
  \begin{eqnarray}
\nonumber
     \IP(N_n(k) \leq k) &=& 1- (1-1/n)\times \ldots (1-k/n)\\ \nonumber
&\leq & 1- (1-k/n)^k\\   \label{eq:45980}
&\sim & k^2/n \qquad {\rm as} \; k^2/n \to 0.
  \end{eqnarray} 
Let $\ep_0=(1/2)(1-1/\E K)>0$. Fix $\ep \in (0,\ep_0)$, and note that  $(1-2\ep) \E K >1$.
It remains to prove the convergence
\begin{equation}
  \label{eq:09221}
  \IP(k \leq \tau_n \leq n \ep, \SGW) =   \sum_{i=k}^{n \ep}   \IP( \tau_n=i,  \SGW)
  \leq \sum_{i=k}^{n \ep}   \IP( S_n(i) \leq 0)
 \to 0,
\end{equation}
where the inequality holds 
by definition of $\tn$ (\ref{eq:tau=min}).
We start to show that there exists a constant $C_\ep >0$, independent of $n$, such that 
$$
\IP\left( N_n(i) < \lceil (1-2\ep)i \rceil \right) \leq \exp \left(-C_\ep i \right)\;,\qquad \forall i \leq n \ep.
$$
Indeed, 
the above probability is equal to 
\begin{equation*}
\IP\left( N_n(i) < \lceil (1-2\ep)i \rceil \right) =  \IP\left( T_{\lceil (1-2\ep)i \rceil,n} > i \right) \leq   \IP \left( \bar T_{(1-2\ep)i}^\ep
> i \right)
\end{equation*}
with $\bar T_{(1-2\ep)i}^\ep$ a sum of $\lceil  (1-2\ep)i \rceil$
i.i.d. geometric r.v.'s with parameter $1-\ep$;
now, the desired estimate follows from Chernov's bound. Next, we note that,
for $z \in (0,1)$, $i \leq n \ep$ and $G(z)=\E z^K$,
\begin{eqnarray*}
  \IP\bigl( S_n(i) \leq 0, N_n(i) \geq \lceil (1-2\ep)i \rceil \bigr) &\leq &
\E \left[ z^{S_n(i)} ;   N_n(i) \geq  \lceil (1-2\ep)i \rceil \right] \\
 &\leq & z^{-i} 
\E\left[ z^{R(\lceil (1-2\ep)i \rceil)} ;   N_n(i) \geq \lceil (1-2\ep)i \rceil \right] \\
 &\leq & z^{-i}
  G(z)^{\lceil (1-2\ep)i \rceil} 
\leq z^{-i} G(z)^{(1-2\ep)i},
\end{eqnarray*}
where we have used $S_{n}(i) = R(N_{n}(i)) - i \geq R(\lceil (1-2\ep)i \rceil )-i$ on the event $\{N_{n}(i) \geq \lceil (1-2\ep)i
\rceil\}$ to pass from the first to the second line, $R$ being given by \eqref{eq:R:2}. Since $(1-2\ep) \E K >1$, we have $r:=z^{-1}G(z)^{(1-2\ep)}<1$ by picking
$z<1$ close enough to 1 and then by expanding $G(z)$ as $G(z) = 1 + \E K (z-1) + o(z-1)$. 
Thus, the left-hand side of (\ref{eq:09221}) 
is bounded by
$$
 \IP(k \leq \tau_n \leq n \ep, \SGW) \leq \sum_{i=k}^{n \ep} 
[r^{i} + \exp (-C_\ep i)] \leq 2(1-r_1)^{-1}  r_1^{k+1},
$$
with $r_1=\max \{r, \exp (-C_\ep)\}<1$. 
Collecting the above estimates in (\ref{eq:tau>}) and taking $k=\lfloor \ln^2 n \rfloor$, 
we conclude that 
$\IP(\tau_n \leq n \ep, \SGW) = {\mathcal O}(n^{-a})$ for all $a \in (0,1)$ (${\mathcal O}(\cdot)$ standing for the Landau notation).
\qed

\noindent  {\bf Proof
 of Theorem \ref{th:LLN}:}
 We start with the proof of (i). We assume first $\E K < \8$. 
Then, we can apply the law of large numbers to the process
\begin{equation} \label{eq:R}
R(m)= \sum_{i=1}^{m} K_i
\end{equation}
in \eqref{eq:R:2}, to show that $\IP$-a.s., $R(nq)/n \to q \E K$ uniformly on $[0,1]$. Recall from \eqref{eq:R:2}
that 
\begin{equation} \label{eq:SRN}
S_n(t)= R(N_n(t))- \lfloor t \rfloor.
\end{equation}
In addition to Lemma \ref{lem:PI},
this shows that, in probability, 
\begin{equation} \label{eq:sommeil}
S_n(ns)/n \longrightarrow (1-e^{-s})\E K -s , \qquad {\rm uniformly\ on \ compacts\ of \ \R_+}
\end{equation}
as $n \to \8$. As a consequence, for any $\delta >0$, 
\begin{equation}
\label{eq:small:event}
\IP \left( \tau_{n} > n (\theta + \delta) \right) \leq
\IP \left( \inf \left[ S_{n}(ns), s \in [0,\theta+\delta] \right] \geq 0 \right) \to 0 \quad {\textrm as} \ n \to \infty,
\end{equation}
since $(1-e^{-s})\E K -s < 0$ for $s > \theta$. 
Now, with $A^\complement$ the complement of $A$, we write 
\begin{eqnarray*}
\IP \left( |\tn -n \theta \;{\bf 1}_{\SGW}|>n \delta \right) &=& 
\IP \left( |\tn -n \theta|>n \delta, \SGW \right) + \IP \left( \tn >n \delta, ({\SGW})^{\complement} \right)
\\
&\leq& \IP\left( \tn < n \ep, \SGW \right) +  \IP \left(    |\tn -n \theta|>n \delta,  \SGW, \tn \geq n \ep \right)\\ 
&& \qquad \qquad \qquad \qquad \qquad 
 \qquad \qquad
+ \IP\left( \tn > n \delta, ({\SGW})^{\complement} \right),
\end{eqnarray*}
where the first term of the right-hand side tends to $0$ (as $n \to \infty$) from Lemma \ref{lem:tau}. 
The last term vanishes because $\tau_n$ is smaller than the extinction time of the Galton-Watson process, which is a.s. finite on the extinction event.
Since $\tn$  
is the first time such that $S_n(\tn)=0$, 
the second term also tends to $0$ by (\ref{eq:sommeil}). Indeed, \eqref{eq:small:event} yields:
\begin{equation*}
\begin{split}
&\limsup_{n \to \infty} \IP \left(    |\tn -n \theta|>n \delta,  \SGW, \tn \geq n \ep \right)
\\
&\quad \leq \limsup_{n \to \infty} \IP \left(  n \ep \leq  \tn \leq n \theta - n \delta,  \SGW, \tn \geq n \ep \right)
\\
&\quad \leq  \limsup_{n \to \infty} \IP \left( \inf \left[ S_{n}(ns), s \in [\varepsilon, \theta - \delta ] \right] =0 \right) 
= 0,
\end{split}
\end{equation*}
since $(1-e^{-s})\E K -s > 0$ for $s \in [\ep, \theta-\delta]$. This ends the proof of the first claim in (i) when $\E K$ is finite. 
The second claim in (i) is then a straightforward consequence of the first one and of (ii) in Lemma \ref{lem:PI}, from which 
$N_{n}(ns)/n \rightarrow 1 - e^{-s}$ in probability, uniformly on compacts of $\R_{+}$.

When $\E K = \8$,  
we consider $K^{(L)}(w)=\min\{K(w),L\}$ for a truncation level $L>0$. By the above proof, we obtain $\lim_{n \to \infty} \tau_n^{(L)}/n = \theta^{(L)} {\mathbf 1}_{{\SGW}^{(L)}}$ with obvious notations. Since
$\tau_n^{(L)} \leq \tau_n$,  $\lim_{L \to \8} \theta^{(L)}=\theta=\8$ and ${\SGW}^{(L)} \nearrow \SGW$ as $L \nearrow \infty$, 
we deduce that $\tau_{n}/n \to \theta$ on $\SGW$ by letting $L$ tend to $\8$.  On the extinction event $({\SGW})^{\complement}$, we obviously have $\tau_{n}/n \to 0$ since $\tau_{n}$ is less than the total number of nodes in the tree. This proves the first claim in (i).
For the second one, we note in the same way that 
$\lim_{n \to \infty} N_{n}(\tau_n^{(L)})/n = p^{(L)} {\mathbf 1}_{{\SGW}^{(L)}}$ and 
$\lim_{L \to \infty} p^{(L)} = 1$. 
Since, $N_n(\tau_n^{(L)})\leq N_n(\tau_n) \leq n$, we deduce that $\lim_{n \to \infty} N_{n}(\tau_n)/n = 1$ on 
$\SGW$. On the extinction event, it obviously holds $N_{n}(\tau_{n})/n \to 0$. 

We now turn to the proof of assertion (ii). By construction, it must hold $N_{n}(\tau_{n}) \leq \tau_{n} +1 \leq 
Z_{\rm tot}^{\rm GW}$, which is a.s. finite when $\E K \leq 1$. By (\ref{eq:45980}), for any $A \geq 0$, 
$\lim_{n \to \8} \IP (N_{n}(\tau_{n}) \leq \tau_{n}, Z_{\rm tot}^{\rm GW} \leq A) = 0$, which proves that, asymptotically,  all the emissions before exhaustion of the capital
are (almost surely) successful on the event $\{Z_{\rm tot}^{\rm GW} \leq A\}$. Therefore, asymptotically, all the emissions before exhaustion are (almost surely) successful, which is to say that, asymptotically with probability 1, every node of the tree receives the information. 
\qed
\subsection{Proof of the Gaussian fluctuations in the case of a light tail}

 {\bf Proof of Theorem \ref{th:flgauss}:}
 By the invariance principle, 
$
(n^{-1/2} \big( R(nq)-nq \E K \big) \cvlaw \bar B(q \si_K^2))_{q \geq 0}
$
in $D(\R_{+})$
with $\bar B$ a Brownian motion. By independence of $(K_i)_i$ and $(\Din)_i$, we have 
from Lemma \ref{lem:PI}:
\begin{equation}
\label{eq:dea0}
n^{-1/2} 
\begin{pmatrix}
( R(nq)-nq \E K )_{q \geq 0}\\  (N_n(ns) -n(1-e^{-s}))_{s \geq 0} 
\end{pmatrix}
 \cvlaw 
\begin{pmatrix}
\bigl[\bar B\big(q \si_K^2\big)\bigr]_{q \geq 0}
\\
\bigl[ B\big( \si_N(s)^2 \big) \bigr]_{s \geq 0}
\end{pmatrix}
\end{equation}
 as $n \to \8$, in $D(\R_{+})^2$ endowed with the product topology generated by the J1 topology (the convergence holding true as well for the uniform topology), where $B$ and $\bar B$ are independent. Then, the convergence
\begin{equation}
	\label{eq:dea1}
	n^{-1/2} \left( R(N_n(ns))-n(1-e^{-s}) \E K \right)_{s \geq 0} \cvlaw \left[ \bar B\left((1-e^{-s}) \si_K^2\right)+
	 B \left((\E K)^2\si_N(s)^2 \right) \right]_{s \geq 0}
\end{equation}
holds in the Skorohod space. 
Actually, we claim that \eqref{eq:dea1} also holds conditionally on $\SGW$, that is under $\IP(\cdot \vert \SGW)$. 
The reason is that the process $n^{-1/2} \big( R(N_n(ns))-n(1-e^{-s}) \E K \big)_{s \geq 0}$ and the event 
$\SGW$ are asymptotically independent, see Lemma \ref{lem:independence:0} right below. 

Now, by Theorem \ref{th:LLN}, 
$$
L_n=\tn -n \theta
$$
is such that $L_n/n \to 0$ as $n \to \8$ in probability under  $\IP( \cdot  | \SGW)$. Using (\ref{eq:dea1}) with $ns=\tn=n\theta +L_n$, we get
\begin{equation}
\label{eq:6:11:1}
n^{-1/2} \left( R(N_n(\tn))-n(1-e^{-\theta -L_n/n}) \E K \right) \cvlaw \bar B\left(p \si_K^2 \right)+
	 B\left((\E K)^2\si_N(\theta)^2 \right)
\end{equation}
under $\IP(\cdot \vert \SGW)$. Now, we emphasize from \eqref{eq:theta} that
\begin{eqnarray}
\tau_{n} - n (1 - e^{- \theta - L_{n}/n} ) \E K  &=& n \theta + L_{n} - n ( 1 - e^{-\theta}) \E K - n e^{-\theta}
( 1 - e^{-L_{n}/n} ) \E K \nonumber
\\
&=& L_{n} -  n e^{-\theta}
( 1 - e^{-L_{n}/n} ) \E K.  
\label{eq:6:11:2}
\end{eqnarray}
Since $\IP(R(N_n(\tn)) = \tn | \SGW) \to 1$ as $n \to \8$, we deduce from \eqref{eq:6:11:1} and \eqref{eq:6:11:2} that, conditionally on survival, 
\begin{equation}
\label{eq:6:11:3}
n^{-1/2} \left( L_n-ne^{-\theta} (1-e^{-L_n/n}) \E K \right) \cvlaw \bar B(p \si_K^2)+
	 B((\E K)^2\si_N(\theta)^2).
\end{equation}
Since
\begin{equation*}
n^{-1/2} L_{n} = \frac{n^{-1/2} ( L_n-ne^{-\theta} (1-e^{-L_n/n}) \E K )}{ 1-e^{-\theta} \E K (1 + 
{\mathcal O} ( L_{n}/n))},
\end{equation*}
we deduce 
from \eqref{eq:6:11:3} that
\begin{equation}
\label{eq:dea2}
n^{-1/2}  L_n  \cvlaw  ( 1-e^{-\theta} \E K)^{-1} [\bar B(p \si_K^2)+
	 B((\E K)^2\si_N(\theta)^2)] \qquad {\rm under} \; \IP(\cdot | \SGW).
\end{equation}
Noting that $1-e^{-\theta} \E K = 1 - (1-p) \E K$, this proves the first claim and the value of $\si_\tau^2$.

To prove the second claim, we note that the left-hand sides in (\ref{eq:dea2})  and in the second line of (\ref{eq:dea0}) jointly converge
as a 2-dimensional vector under the conditional law $\IP(\cdot \vert \SGW)$. Therefore, 
$$
n^{-1/2} \Big(  N_n(\tn) -n(1-e^{-\theta}) - e^{-\theta} L_n\Big)
 \cvlaw 
 B\big( \si_N(\theta)^2 \big) \qquad {\rm under} \; \IP(\cdot | \SGW),
$$
and also
$$
n^{-1/2} \Big(  N_n(\tn) -n(1-e^{-\theta})\Big)  \cvlaw ( 1-e^{-\theta} \E K) ^{-1}
\Big[ e^{-\theta} \bar B(p \si_K^2)+
	 B( \si_N(\theta)^2)\Big],
$$
conditionally on $\SGW$. 
This completes the proof.
 \qed

We finally prove
\begin{lem}
\label{lem:independence:0}
Under the assumptions of Theorem \ref{th:flgauss}, 
the process 
$(n^{-1/2}( R(nq)-nq \E K ))_{q \geq 0}$ has the same limits in law under 
$\IP(\cdot \vert (\SGW)^{\complement})$  and $\IP$. 
\end{lem}
 {\bf Proof:}. It is sufficient to prove that, when $\IP(\SGW) < 1$, the process 
$(n^{-1/2}( R(nq)-nq \E K ))_{q \geq 0}$ has the same limits in law under 
$\IP(\cdot \vert (\SGW)^{\complement})$  and $\IP$. 
From \eqref{def:ki}, we know that, 
on the event $\{ Z_{\rm tot}^{\rm GW} \leq A\}$,
for $A \geq 0$, 
 the variable $K_{i}$ in the definition of $R$ (see \eqref{eq:R:2}) coincides with $\bar{K}_{i}$ for 
$i > A$, so that the process $\sum_{j = A}^{\lfloor n q \rfloor} \bar{K}_{j} = R(nq) - \sum_{j=1}^A K_{j}$ is independent of 
the Galton-Watson tree. Therefore, the event $\{Z_{\rm tot}^{\rm GW} \leq A\}$ and the process
$(n^{-1/2}( R(nq)-nq \E K ))_{q \geq 0}$ are asymptotically independent, which is to say that $(\SGW)^{\complement}$
and $( n^{-1/2}(R(nq)-nq \E K) )_{q \geq 0}$ are asymptotically independent. Therefore, 
the process 
$(n^{-1/2} (R(nq)-nq \E K ))_{q \geq 0}$ has the same limits in law under 
$\IP(\cdot \vert (\SGW)^{\complement})$  and $\IP$.
\qed
\subsection{Proof of the fluctuations in the case of a heavy tail}

 {\bf Proof of Proposition \ref{th:FT+CT}:} Clearly, $\E K=\infty$. By Theorem \ref{th:LLN}, for any $\varepsilon >0$, as $n \rightarrow \infty$,
$$\IP(N_{n}(\tau_{n})/n \geq 1-\varepsilon, \SGW) \longrightarrow 
\IP(\SGW).$$ 
Moreover, 
$\{N_{n}(\tau_{n}) = n\} \supset \{N_{n}(\tau_{n}) \geq n/2, R(\lfloor n/2 \rfloor) \geq T_{n,n}\}$ since 
$R(\lfloor n/2 \rfloor) \geq T_{n,n} \Rightarrow R(k) \geq k$ for any $k \in \{\lfloor n/2 \rfloor, \dots, T_{n,n}\}$. 
Thus, for any $\beta \in (1,1/\alpha)$,
\begin{equation*}
\begin{split}
\liminf_{n\rightarrow \infty} \IP\bigl( N_{n}(\tau_{n}) = n,\SGW\bigr) 
\geq \liminf_{n\rightarrow \infty}  \IP \left( R(n/2) \geq n^{\beta},
T_{n,n} \leq n^{\beta},\SGW\right). 
 \end{split}
 \end{equation*} 
 By Markov inequality, $\IP ( T_{n,n} > n^{\beta} ) \rightarrow 0$ as $n \rightarrow \infty$, since 
 $\E T_{n,n} = \sum_{i=1}^{n-1} n/(n-i) \sim_{n \to \8} n \ln (n)$.
Moreover, $\IP(K_{1} \geq n^{\beta}) \geq (c/2) n^{-\alpha \beta}$ for $n$ large enough, so that, for $n$ large, 
\begin{equation*}
\IP \biggl( \bigcap_{i=1}^{\lfloor n/2 \rfloor} \{K_{i} < n^{\beta}\} \biggr) \leq 
\bigl( 1 - \frac{c}{2 n^{\alpha \beta}} \bigr)^{\lfloor n/2 \rfloor} \sim \exp\bigl( - \frac{c}{4} n^{1-\alpha \beta} \bigr)
\longrightarrow 0 \qquad \textrm{as} \ n \rightarrow \8.
\end{equation*}
Therefore, $\IP ( \sum_{i=1}^{\lfloor n/2 \rfloor } K_{i} \geq n^{\beta}) \rightarrow 1$. 
We deduce that 
$\liminf_{n\rightarrow \infty} \IP ( N_{n}(\tau_{n}) = n,\SGW ) = \IP(\SGW)$. 
Finally, on $(\SGW)^{\complement}$, we have  $N_{n}(\tau_{n})/n \rightarrow 0$ in probability, so that
 $\IP ( N_{n}(\tau_{n})=n $,$(\SGW)^{\complement}) \rightarrow 0$. \qed

\noindent
 {\bf Proof of Theorem \ref{th:critique}:}
By a celebrated result of \cite{ErRe61} (see \cite{Dur95} pp. 130-132, for a short account), 
\begin{equation}
\label{eq:erdosrenyi}
G^{(n)}=n^{-1}\big(T_{n,n} - n \ln n\big) \cvlaw G \;,
\end{equation}�
where the variable $G$ has a Gumbel distribution, 
\begin{equation}
\label{eq:gumbel}
\IP(G \leq x)=e^{-e^{-x}}, x\in \R. 
\end{equation}� 
On the other hand, from the tail assumption for $K$,
\begin{equation}
\label{eq:cauchy}
 {\mathcal S}_c^{(n)}=  n^{-1}\big(R(n-1)-cn\ln n\big) \cvlaw {\mathcal S}_c\,,
\end{equation}
where $ {\mathcal S}_c$ is a totally asymmetric, stable law with index 1 (Cauchy law), depending upon the parameter $c$.
We define  ${\mathcal S}= {\mathcal S}_1$. Precisely, the law of  ${\mathcal S}$ is given by
$$
\E e^{iu{\mathcal S}} = \exp \biggl\{ \int_0^1 (e^{iux}-1-iux) x^{-2}dx
+ \int_1^\8 (e^{iux}-1) x^{-2}dx + iuc_0 \biggr\}, \quad u \in \R,
$$
where $i^2=-1$ and $c_0 \in \R$ is defined by $c_0=\lim_{n \to \8} (\E[ K; K \leq n] - c \ln n)$.
Recall from Theorem \ref{th:LLN} that $\IP(\Trans_n \cap \SGW)-\IP(\Trans_n) \to 0$.
We have
\begin{eqnarray*}
\Trans_n \cap \SGW &=& \{\tn \geq T_{n,n}\} \cap \SGW\\
&\subset &  \{ R(n-1) \geq T_{n,n}\} \cap \SGW\\
&\subset & 
 \{ G^{(n)} -  {\mathcal S}_c^{(n)}
\leq (c-1)\ln n
\} \cap \SGW,
\end{eqnarray*}
where we have used $\tau_{n} \geq T_{n,n} \Rightarrow 
R(n-1) = R(N_{n}(T_{n,n}-1)) > T_{n,n}-1$ to pass from the first to the second line. 
The random vector $( G^{(n)},  {\mathcal S}_c^{(n)})$ converges in law to a couple  $( G,  {\mathcal S}_c)$ with independent components (independence follows from the independence of the processes $N_{n}$ and $R$).   
If $c<1$, then $\IP(G^{(n)} -  {\mathcal S}_c^{(n)}
\leq (c-1)\ln n) \rightarrow 0$ so that $\IP(\Trans_n) \rightarrow 0$. 
If $c>1$, it obviously holds $\limsup_{n \rightarrow \8}
\IP(\Trans_n) \leq \IP(\SGW)$. 
To tackle the case when $c=1$, we observe
from Lemma \ref{lem:independence} below
 that 
$\SGW$ and ${\mathcal S}_1^{(n)}$ are asymptotically independent 
as $n \to \8$, so that the random vector $( G^{(n)},  {\mathcal S}_1^{(n)},{\mathbf 1}_{\SGW})$ converges in law to a triple 
$(G,{\mathcal S},{\mathbf 1}_{A})$ with independent components, where $\IP(A) = \IP(\SGW)$. 
As the random variable $G-{\mathcal S}$ has a continuous cumulative distribution function, we get:
\begin{equation*}
\begin{split}
\limsup_{n \to \8} \IP(\Trans_n) &\leq
\limsup_{n \to \8} \IP \bigl(  \{ G^{(n)} -  {\mathcal S}_{1}^{(n)}
\leq 0
\} \cap \SGW \bigr)
\\
&=
\IP(\SGW)  
\IP( G \leq {\mathcal S}) = \IP(\SGW) \E ( \exp(- e^{-{\mathcal S}})),
\end{split}
\end{equation*}
the last equality following from \eqref{eq:gumbel}. 

We now turn to the reverse bound in the case when $c \geq 1$ (when $c<1$, the proof is over). Consider $\varepsilon \in (0, 1)$ and a positive sequence $(\delta_{n})_{n \geq 1}$ such that 
$$\delta_{n} \ln n \rightarrow 0, \qquad n \delta_n \to \8,$$
as $n \rightarrow \8$.  
We write
 \begin{eqnarray} \nonumber
\Trans_n \cap \SGW \supset 
 \{ R(n(1-\delta_n))>T_{n,n}, \SGW \}
&\cap &
 \{N_n(\tn) \geq n(1-\varepsilon), \SGW \}
\\ & \cap& \{ R(n(1-\varepsilon))>T_{n(1-\delta_n),n}\}, \nonumber
\end{eqnarray}
the right-hand side being denoted $A \cap B \cap C$. 
Indeed, on $B \cap C$, $\tau_{n} < T_{n,n} \Rightarrow \tau_{n} = R(N_{n}(\tn)) > T_{n(1-\delta_{n}),n}$ so that 
$\tn > T_{n(1-\delta_{n}),n}$. Then, by the same argument, it must hold $\tn \geq T_{n,n}$ on $A \cap B \cap C$ as otherwise  
$\tn$ would be equal to $R(N_{n}(\tn)) \geq R(N_{n}(T_{n(1-\delta_{n}),n})) = R(n(1-\delta_{n})) > T_{n,n}$, yielding to a contradiction. Below we will  use the estimate $\IP(A \cap B \cap C) \geq \IP(A)-\IP(B^\complement \cap A) -\IP(C^\complement)$.

As above, we write $ \{ R(n(1-\delta_n))>T_{n,n}, \SGW \}$ as
$\{ n^{-1} (R(n(1-\delta_n)) - c n \ln n) - G^{(n)} > - (c-1) \ln n, \SGW\}$. 
Lemma \ref{lem:independence} says that $n^{-1} (R(n(1-\delta_n)) - c n \ln n) - G^{(n)}$ converges in law towards 
${\mathcal S}_{c}-G$ as $n \rightarrow \8$. Therefore, when $c>1$, we get 
\begin{equation*}
\lim_{n \rightarrow \8} \IP \bigl( R(n(1-\delta_n))>T_{n,n}, \SGW \bigr)
= \IP (\SGW). 
\end{equation*} 
When $c=1$, we make use of Lemma \ref{lem:independence} again. By asymptotic independence, we get as in the proof of the upper bound:
\begin{eqnarray*}
\lim_{n \to \8} \IP \left( R(n(1\!-\!\delta_n))\!>\!T_{n,n}, \SGW \right)
&=&\lim_{n \to \8} \IP\Bigl( (\frac{R(n(1\!-\!\delta_n))}{n}- \ln n) - G^{(n)} >0, \SGW \Bigr)
\\
&= & \IP(\SGW) 
\IP( G < {\mathcal S}) = \IP(\SGW) \E ( \exp(- e^{-{\mathcal S}})).
\end{eqnarray*}
 By Theorem \ref{th:LLN} with $\E K= \8$, 
$\IP(N_{n}(\tau_{n}) \leq n(1-\varepsilon), \SGW) \rightarrow 
0$ as $n \rightarrow \infty$. It remains to show that 
$$
\IP\left(  R(n(1-\varepsilon))>T_{n(1-\delta_n),n} \right) \to 1 \qquad \textrm{as} \quad n \to \8.
$$
We have
\begin{eqnarray*}
\E \left( T_{n(1-\delta_n),n} \right) &=& \sum_{i=1}^{\lfloor n(1-\delta_n) \rfloor -1} \frac{n}{n-i} \sim n \ln 1/\delta_n, \\
{\rm Var} \left( T_{n(1-\delta_n),n} \right) &\leq & 
\sum_{i=1}^{\lfloor n(1-\delta_n) \rfloor} \frac{n^2}{(n-i)^2} \sim n/\delta_{n} = o(n^2),
\end{eqnarray*}
as $n \delta_{n} \to \8$. Tchebyshev inequality implies that $T_{n(1-\delta_n),n}= {o}(n \ln n)$ for $n \delta_n \to \8$.
On the other hand, $ R(n(1-\varepsilon))$ is of order $n \ln n$, which ends the proof. \qed
\vspace{5pt}

We finally prove
\begin{lem}
\label{lem:independence}
Under the assumptions of Theorem \ref{th:critique}, 
given a sequence of positive reals $(\delta_{n})_{n \geq 1}$ such that $\delta_{n} \ln n \rightarrow 0$ as $n \rightarrow \8$, the sequence 
$(n^{-1}R (n(1- \delta_{n})) - c \ln (n))_{n \geq 1}$ converges in law towards 
${\mathcal S}_{c}$. Moreover, as $n \rightarrow \infty$, the event $\SGW$ and the variable $n^{-1} R(n(1-\delta_{n})) - c \ln n$ become independent.
\end{lem}
 {\bf Proof:}. In order to prove that the sequence $(n^{-1}R (n(1- \delta_{n})) - c \ln (n))_{n \geq 1}$ converges in law towards 
${\mathcal S}_{c}$, it is sufficient to check that $(n^{-1}[R(n-1) - R(n(1-\delta_{n}))])_{n \geq 1}$ converges towards $0$ in probability. For any $\varepsilon >0$ and $A>0$, we have:
\begin{equation*}
\begin{split}
\IP \bigl( R(n-1) - R(n(1-\delta_{n})) \geq n \varepsilon \bigr) 
&\leq \IP \bigl( \exists \ell = \lfloor n (1-\delta_{n}) \rfloor + 1, \dots, n-1 : K_{\ell} \geq A \bigr)
\\
&\hspace{15pt} +
\IP \biggl( \sum_{i=	\lceil n(1-\delta_{n}) \rceil}^n K_{i} {\mathbf 1}_{\{K_{i} \leq A\}} \geq n \varepsilon
\biggr) 
\\
&\leq \IP \bigl( \exists \ell = \lfloor n (1-\delta_{n}) \rfloor + 1, \dots, n-1 : K_{\ell} \geq A \bigr) 
\\
&\hspace{15pt} + \varepsilon^{-1} \delta_{n} {\mathbb E} \bigl[ K {\mathbf 1}_{\{K \leq A\}} \bigr]
\\
&\leq 1 - \bigl( 1 - \IP(K \geq A) \bigr)^{n \delta_{n}} + c'  \varepsilon^{-1} \delta_{n} \ln (A), 
\end{split}
\end{equation*}
for some constant $c'$ independent of $n$ and $A$. For $A = n$, we have $1 - \IP(K \geq n) \sim 1 - c/n$ by assumption 
and thus  
$( 1 - \IP(K \geq n) )^{n \delta_{n}} \sim \exp( - c \delta_{n})$, which tends to $1$ as $n \rightarrow \8$.  
As $\delta_{n} \ln n \rightarrow 0$ as $n \rightarrow \8$, we deduce that $(n^{-1}[R(n-1) - R(n(1-\delta_{n}))])_{n \geq 1}$ indeed converges towards $0$ in probability.

The asymptotic independence is proved as in the proof of Lemma \ref{lem:independence:0}.
\qed
\vspace{5pt}


\section{Proofs of large deviations} \label{sec:prld}


Before starting the proofs we recall a few facts about the LDP for the coupon collector.
As $r$ is increased from 0 to $t \wedge 1$, the function $r \mapsto \rho(r,t)$ in (\ref{eq:rhodef}) decreases from $+\8$ to  zero when $t\leq 1$ or to a positive value otherwise.
The value $\rho=1$ corresponds to the typical case  $r=1-e^{-t}$; and $\rho >1$
(resp. $\rho < 1$) 
to large deviations with $N_n(nt)$ much smaller (resp. larger) than $ \E[N_n(nt)]$ (cf Appendix A.2.1 of \cite{DNW}).
The LDP holds at the process level for $(N(nt); t \geq 0)$, and the overwhelming  contribution to the probability
of the event $\{N_n(nt) \in n(r-\ep, r+\ep)\}$ occurs in the neighborhood of the optimal path $N_n(ns)=ng(s)$ with
\begin{equation} \label{eq:02}
g(s)= \frac{1-e^{-s\rho}}{\rho},\qquad s\in [0,t],
\end{equation}
and  $\rho=\rho(r,t)$ as before (see also Lemma \ref{lem:uniqueness} below for the uniqueness of the optimal path). 
\medskip

 {\bf Proof of Lemma \ref{lem:continuity}:}
(i) Lower semi-continuity is easily checked. Boundedness of the level sets follows from the lower bound
\begin{equation*}
r I\bigl(\frac{t}{r} \bigr) \geq r \bigl( a_{0} \frac{t}{r} - \ln \E \exp(a_{0} K) \bigr) 
\geq a_{0} t - \ln \E \exp(a_{0} K) \bigr),
\end{equation*}
which implies for a well-chosen value of $a_{0}$ that $t$ must be bounded when ${\mathcal F}(r,t)$ is bounded; boundedness of 
$r$ easily follows since $r \leq t$ when ${\mathcal F}(r,t)$ is finite. 

To determine the domain of ${\mathcal F}$, we recall the expression of  $ {\rm Dom}(I)$ in terms of $k_*, k^*$, and  write
$${\rm Dom}({\mathcal F})=\{t/r \in {\rm Dom}(I), r/t \in ]0,1]\} \cup \{(0,0)\}.
$$
Continuity on the first set in the above union follows from that of $I$ on its domain. 
If $k^*=\8$, $\cF$ is unbounded in any neighborhood of the origin, so it is not continuous at this point.
If $k^*<\8$, then $I$ is bounded on its domain and continuity at the origin of $\cF$ easily follows.

(ii) trivially holds. 
(iii)
Finally, roots of ${\mathcal F}$ must satisfy $r=0$ or $t/r=\E K$, for $0 < r \leq t \wedge 1$, as $\E K$ is the only zero of $I$. 
When $r=0$, $t$ must be zero as well. When $0< r \leq t \wedge 1$, the condition $t/r = \E K$ requires $\E K$ to be strictly greater than 1. (When $\E K=1$, the condition $t=r>0$ implies $t \leq 1$ and $J_{t}(r)=J_{t}(t)=t >0$.) When $\E K >1$, $t/r=\E K$ 
implies $t \rho= \theta$ (compare with   \eqref{eq:theta}) and $J_{t}(r)=0$ implies $r=1 - e^{-t}$, that is $\rho=1$ and thus 
$t= \theta$ and $r=p$ (see \eqref{eq:p(theta)}). Conversely, it is well-checked that $0$ is a root of ${\mathcal F}$ and that 
$(p,\theta)$ is a root as well when $\E K > 1$. 
 \qed

\vspace{5pt}

\noindent  {\bf Proof of Theorem \ref{th:ldp}:}

{\bf Upper bound}. 
We start with the case when $r,t$ satisfy $0 < r < t \wedge 1$. 
Given $\varepsilon,\delta >0$, we are to prove the local upper bound
\begin{equation}
\label{eq:ld5}
\limsup_{\eps, \delta \to 0^+} 
\limsup_{n \to \8} \un \ln \IP \bigl( N_n(\tau_n)/n \in [r-\eps,r+\eps], \tau_n/n \in [t-\delta,t+\delta] \bigr) \leq
- \cF(r,t).
\end{equation}
We first tackle the case $r<1-e^{-t}, t/r<\E K$, which is relevant only when $\E K >1$. 
Without any loss of generality, we can assume that $\eps$ and $\delta$ satisfy 
\begin{equation}
\label{eq:ld2}
r+\eps<1-e^{-(t-\delta)},\qquad (t+\delta)/(r-\eps)<\E K.
\end{equation}
We then define the events
\begin{equation}
\label{eq:ld3} 
A=\{N_n(n(t-\delta))\leq n(r+\eps)\},\qquad B=\{R(n(r-\eps))\leq n(t+\delta)\}.
\end{equation}
Since $N_n, R$ are nondecreasing,
$$
\bigl\{N_n(\tau_n)\in [n(r-\eps),n(r+\eps)], \tau_n \in [n(t-\delta),n(t+\delta)]\bigr\} \subset A \cap B,
$$
and by independence of the processes $N_{n}$ and $R$, we get
\begin{equation}
\label{eq:ld4}
\IP \bigl( N_n(\tau_n)/n \in [r-\eps,r+\eps], \tau_n/n \in [t-\delta,t+\delta] \bigr) \leq \IP(A) \IP( B).
\end{equation}� 
From the LDP's (\ref{eq:ld0}) and (\ref{eq:ld1}), and by (\ref{eq:ld2}), we have
\begin{equation*}
\begin{split}
&\limsup_{n \to \8} \un \ln \IP \bigl( N_n(\tau_n)/n \in [r-\eps,r+\eps], \tau_n/n \in [t-\delta,t+\delta] \bigr) 
\\
&\hspace{15pt}\leq
-  J_{t-\delta}(r+\eps) - (r-\eps) I((t+\delta)/ (r-\eps)). 
\end{split}
\end{equation*}
By lower semi-continuity  of $I$ on the whole $\R$ and by continuity of $(r,t) \mapsto J_t(r)$ 
at the prescribed value of $(r,t)$, we obtain
\eqref{eq:ld5}.
With similar arguments, one easily obtain the same result in all other cases of $0<r<t\wedge 1$.

We are now left with proving the local upper bound on the boundary. For $(r,t)=(0,0)$, there is nothing to prove since the rate function is zero. The case $r=0<t$ is simple since it is enough to bound
$$
\IP \bigl( N_n(\tau_n)/n \leq\eps, \tau_n/n \in [t-\delta,t+\delta] \bigr)
\leq \IP \bigl( N_n(n(t-\delta))/n \leq\eps \bigr),
$$
to get a rate of decay $J_{t-\delta}(\eps)$ which tends to $\8$ as both $\delta$ and $\eps$ vanish. 
The cases $r=1<t$ and $r=t \leq 1$  use similar arguments as above in the general case. 

We have proved the upper bound for compact sets. To extend it to closed sets, it is enough to show that
\begin{equation}
\label{eq:ld6} 
\lim_{t \to \8} \limsup_{n \to \8} \un \ln \IP\bigl(  \tau_n/n \geq t\bigr) =-\8.
\end{equation}
\medskip
For this, we observe that
$$
\{\tn \geq nt\} \subset \{R(N_n(nt)) \geq \lfloor nt \rfloor \}  \subset \{R(n)\geq \lfloor nt \rfloor \}.
$$
For $t > \E K$, the probability of the last event can be estimated by Cram\'er's bound, yielding
$$
\IP( \tn \geq nt ) \leq \exp (-n I(t)).
$$
Since $I(t)$ tends to $\8$ as $t \to \8$, we obtain (\ref{eq:ld6}).
\bigskip

{\bf Lower bound}. 
The lower bound is subtle. It is sufficient to show that, for all $r,t$ with $\cF(r,t) < \infty$ and all $\eta >0$,
\begin{equation}
\label{eq:lower:bound:0}
\liminf_{n \to \8} \un \ln
\IP \bigl(   N_n(\tau_n)/n \in [r-\eta, r+\eta], \tau_n/n \in [t-\eta, t+\eta] \bigr) \geq 
-\cF(r,t).
\end{equation}

We start with the general case when $r,t$ satisfy $0<r<t\wedge 1$ and $I(t/r) < \infty$. (If $I(t/r)=\infty$, the lower bound above is trivial.)

\textit{First Step.} The proof holds in several steps. The first one is to bound from below the left-hand side above
by the probability of an event depending in a separate way on the dynamics of the coupon collector on the one hand and on the capitals of the servers and the Galton-Watson tree on the other hand. We are thus given $0 < \eta' < \eta < r$ and $0 <  \eps' < 2 \eps$ such that 
$r + \eta + \eps < 1 \wedge t$, $\varepsilon < r\eta/(4t) < \eta/4$, $\eta' < (1- r/t) \eta$ 
 and $\eta' (1+t/r) < \varepsilon'$.
For an integer $m \in [n \eps',n (\eps-\eps')]$,
we define the events
\begin{equation*}
\begin{split}
&A=\{N_n(nt+ n\eta)\leq nr + n\eps \}, 
\\
&A^{\prime}=\{N_n(n\eps)= \lfloor n\eps \rfloor +1\}, 
\\
&A^{\prime \prime}=\{N_n(\ell)\geq \ell r/t + n \varepsilon ( 1 - r/t) - n \eta' ; n\eps \leq \ell \leq nt-n\eta\}
\cap \{ N_{n}( nt - n \eta) \leq n r + n \eps \},
\\
&B=\{ R(nr + n\eps)-R(m) \leq nt+n\eta/2\}, 
\\
&B^{\prime \prime}=\{ R(\ell)-R(m) \geq (\ell-m)t/r - n \eta';  m 
 \leq \ell \leq nr + n\eps\}.
\end{split}
\end{equation*}
Recall that $N_n(s)=N_n(\lfloor s \rfloor)$ and similarly for $R$. We also define
\begin{equation}\label{def:AAA2}
B'=\{ R(m) \in [n\eps,2 n \eps],
Z_{\rm tot}^{\rm GW} \geq n\eps\}.
\end{equation}
On $A'$, it holds $T_{i+1,n}=i (i \leq n\eps)$, so that the emission process cannot stop before $n \varepsilon$ unless the total size of the Galton-Watson tree is strictly less than $n \varepsilon$. Therefore, $A'\cap B' \subset \{ \tau_n \geq n\eps\}$. Further,
$A'\cap A^{\prime \prime} \cap B'  \cap B^{\prime \prime} \subset \{ \tau_n > n (t-\eta)\}$. 
Indeed, on $A'\cap A^{\prime \prime} \cap B'  \cap B^{\prime \prime}$, for $n \varepsilon \leq \ell \leq n t -  n\eta$, 
\begin{equation}
\label{eq:ld:200}
\begin{split}
&R \bigl( N_{n}(\ell) \bigr) 
\\
&\geq R \bigl( 
\ell r/t + n \varepsilon ( 1 - r/t) - n \eta' \bigr) 
\\
&\geq 
\bigl( \ell r/t + n \varepsilon ( 1 - r/t) - n \eta' - m  \bigr) t/r + R(m) -  n \eta'
\\
&\geq \ell + (n \varepsilon - m)(t/r- 1) + R(m) - m  - n \eta'(1+t/r) > \ell + n \varepsilon' -  n \eta'(1+t/r) > \ell, 
\end{split}
\end{equation}
where, to pass from the second to the third line, we use the fact that, for the prescribed values of $\ell$, 
$\ell r/t + n \varepsilon ( 1 - r/t) - n \eta' \in [ m, n r + n \eps]$.

Moreover, on $A \cap B \cap B'$, 
\begin{equation*}
R \bigl( N_{n}(n t + n \eta) \bigr) \leq R \bigl( nr + n\eps  \bigr) \leq 
n t + n \eta/2 + R(m) \leq nt + n \eta/2 + 2 n \eps < nt + n \eta,
\end{equation*}
so that, $\tau_{n} \leq n (t+ \eta)$
on $A \cap B \cap B'$. Therefore,
\begin{equation}
\label{eq:inclusion}
\big( A \cap A' \cap A^{\prime \prime} \cap B \cap B' \cap B^{\prime \prime} \big) \subset
\{  N_n(\tau_n) \in [n(r-\eta), n(r+\eta)], \tau_n \in [n(t-\eta), n(t+\eta)]\},
\end{equation}
since, on $A \cap A' \cap A^{\prime \prime} \cap B \cap B' \cap B^{\prime \prime}$,
\begin{equation}
\label{eq:ldp:300}
\begin{split}
&N_{n}(\tau_{n}) \geq N_{n}(n(t-\eta)) \geq n (t-\eta) r /t + n \varepsilon (1 - r /t) - n \eta' 
\\
&\hspace{109pt} \geq n r - n \eta r/t - n \eta' \geq nr - n \eta,
\\
&N_{n}(\tau_{n}) \leq N_{n}(n(t+\eta)) \leq nr + n\eps \leq n r + n \eta. 
\end{split}
\end{equation}
By independence of the processes $N_{n}$ and $R$, we get
\begin{equation}
\label{eq:ld11}
\IP(   N_n(\tau_n)/n \in [r-\eta, r+\eta], \tau_n/n \in [t-\eta, t+\eta]) \geq 
\IP\big( A \cap A' \cap A^{\prime \prime}\big) \IP \big(  B \cap B' \cap B^{\prime \prime} \big). 
\end{equation}
\textit{Second Step.}
We estimate the first factor in the right-hand side by using Markov's property:
\begin{equation}
\label{eq:ld12}
\begin{split}
\IP\big( A \cap A' \cap A^{\prime \prime}\big) &\geq \inf_{i \in \hat{I}}\IP \big( A \vert \hat{A}_{i}
 \bigr)  \IP \big( A'' \vert A' \big) \IP \big( A' \big),
\\
\textrm{with} \ &\hat{A}_{i} = \bigl\{ N_{n}(nt - n \eta) = i \bigr\},
\\
&\ \hat{I} = \{ \lceil n(r-\eta r/t) + n \varepsilon (1 - r/t) - n \eta' \rceil, \dots, \lfloor nr + n \varepsilon
\rfloor \}.
\end{split}
\end{equation}

Clearly, $ \IP( A')$ is bounded from below by the probability for the binomial distribution with parameters $\lfloor n \eps \rfloor$ and $\eps$
to be equal to zero. We get 
\begin{equation}
\label{eq:ld:101}
\IP \big( A' \big) \geq e^{-n \eps |\ln (1-\eps)|}.
\end{equation}
Similarly, $\inf_{i \in \hat{I}}\IP(A \vert \hat{A}_{i})$ is bounded from below by the probability for the binomial distribution with parameters $\lfloor 2 n \eta \rfloor +1$
and $1 - r + \eta r /t - \eps ( 1 - r/t) + \eta' < 1 - r + \eta$ to be equal to $0$. Therefore,
\begin{equation}
\label{eq:ld:102}
\inf_{i \in \hat{I}}\IP \big( A \vert  \hat{A}_{i} \bigr) \geq e^{- (2 n \eta +1) \vert \ln(r-\eta) \vert }.
\end{equation}
We now turn to $\IP(A'' \vert A')$. We observe that the conditional probability $\IP\big( \cdot \vert A' \big)$ is (up to a shift in time) the law of the coupon collector with $n$ images when starting from $n \eps$ different images, already collected at the initial time.
Since $n (r-\eta r/t) + n \eps ( 1 - r/t) - n \eta' \leq n r + n (\eps - \eta r /t) - n \eta' \leq nr - n \eta'$, we deduce from the large deviation lower bound in \cite{DNW}:
\begin{equation}
\label{eq:ld20}
\begin{split}
&\IP\big( \hat{A}
 \vert A' \bigr) \geq  \exp \{-n J^\eps_{t-(\eta+\eps)}(r)+o(n) \}
\\
&\qquad \textrm{with} \  \hat{A} = \bigl\{
n (r- \eta r/t) + n \eps ( 1-r/t) - n \eta' \leq 
N_{n}(nt -n \eta) \leq nr + n\eps \bigr\} = \bigcup_{i \in \hat{I}} \hat{A}_{i}, 
\end{split}
\end{equation}
the rate function $J^\eps_{t-(\eps+\eta)}(r)$ being given by Theorem 2.7 in \cite{DNW} (with $I=0$, $\alpha_{0} = 1- \varepsilon$ and $\omega_{0} = 1-r$ therein).  Precisely, (\ref{eq:ld20}) follows from the LDP for the ``time-shifted'' variable $N_{n}(n(t-(\eta+\eps)))$ with $N_{n}(0)=\lfloor n \eps \rfloor +1$ as initial condition. Using the same notations as in \cite{DNW}, $J^{\eps}_{t}(r)$ may be expressed as the relative entropy:
\begin{equation*}
\begin{split}
J^{\eps}_{t}(r) &= (1- \varepsilon) \bigl[ \frac{1-r}{1- \varepsilon} \ln \bigl( \frac{1-r}{(1-\varepsilon) e^{-t}} \bigr) 
\bigr]
+ (1 - \varepsilon)  \Ceps \sum_{j \geq 1} {\mathcal P}_{j}(\rho^{\eps} t) \ln \bigl[ \Ceps \frac{{\mathcal P}_{j}(\rho^{\eps} t) }{{\mathcal P}_{j}(t)}\bigr]
\\
&\hspace{15pt} + \varepsilon   \sum_{j \geq 0} {\mathcal P}_{j}(\rho^{\eps} t) \ln \bigl[  \frac{{\mathcal P}_{j}(\rho^{\eps} t) }{{\mathcal P}_{j}( t)}\bigr],
\end{split}
\end{equation*}
where 
${\mathcal P}_{j}(\alpha)$ is the $j$th weight of the Poisson distribution of parameter $\alpha$, 
$\rho^{\varepsilon}$ is the unique root of
\begin{equation}
\label{eq:rho:eps}
\frac{1-r}{1-\varepsilon} + \frac{1- \varepsilon \rhops}{(1-\varepsilon)\rhops} (1 - e^{- \rhops t}) = 1,
\end{equation}
and $C^{\varepsilon}$ reads:
\begin{equation*}
\Ceps = \frac{ 1 - \varepsilon \rhops}{(1- \varepsilon) \rhops}. 
\end{equation*}
It is plain to check that, as $\varepsilon$ tends to $0$, $\rhops$ converges towards $\rho=\rho(r,t)$ and $C^{\varepsilon}$ 
towards $1/\rho$. Moreover, standard computations yield
\begin{equation*}
\begin{split}
& \sum_{j \geq 0} {\mathcal P}_{j}(\rho^{\eps} t) \ln \bigl[  \frac{{\mathcal P}_{j}(\rho^{\eps} t) }{{\mathcal P}_{j}( t)}\bigr] = -(\rhops - 1) t + \rhops \ln (\rhops) t ,
\\
&\sum_{j \geq 1} {\mathcal P}_{j}(\rho^{\eps} t) \ln \bigl[  \frac{{\mathcal P}_{j}(\rho^{\eps} t) }{{\mathcal P}_{j}( t)}\bigr] = -(\rhops - 1) t + \rhops \ln (\rhops) t + (\rhops - 1) t e^{-\rhops t},  
\end{split}
\end{equation*}
so that, by the definition of $\rho(r,t)$, 
\begin{equation*}
\begin{split}
\lim_{\varepsilon \rightarrow 0} J^{\eps}_{t}(r) &= (1-r) \ln(1-r) + (1-r) t - \frac{\ln \rho}{\rho} \bigl(1 - e^{-\rho t}\bigr)
- \frac{\rho - 1}{\rho} t + t \ln  \rho   + \frac{\rho - 1}{\rho} t e^{-\rho t}
\\
&= J_{t}(r). 
\end{split}
\end{equation*}
Similarly, it holds 
\begin{equation*}
J^\eps_{t-(\eta+\eps)}(r)  \longrightarrow J_{t-\eta}(r) , \qquad {\rm as} \; \eps \to 0.
\end{equation*}
Therefore, \eqref{eq:ld20} may be expressed as 
\begin{equation}
\label{eq:ld105}
\IP\big( \hat{A} \vert A' \bigr) 
\geq  \exp \{-n [J_{t}(r) + \delta(\eta) + \delta'(\eps;\eta)]
+o(n) \},
\end{equation}
where $\delta(\eta)$ is a generic term that tends to $0$ with $\eta$ and $\delta'(\varepsilon;\eta)$ is a generic term that tends to $0$ with $\eps$ when $\eta >0$ is given. 
We then claim that $\IP(A'' \vert A')$ satisfies the same lower bound, that is
\begin{equation}
\label{eq:ld106}
\IP\big( A'' \vert A' \bigr) 
\geq  \exp \{-n [J_{t}(r) + \delta(\eta) + \delta'(\eps;\eta)]
+o(n) \}.
\end{equation}
Basically, it comes from the fact that the optimal path explaining the LD of the random variable
$N_n(n(t-\eta))$ in the neighborhood of $r$ given the initial condition in $A'$ is lying above the linear constraint in $A''$. As already explained, the right optimal path to consider is the one when the collector has already collected $
\lfloor n \varepsilon \rfloor +1$ coupons at (rescaled) time $\varepsilon$. We denote it by $g^{\varepsilon} : [\varepsilon,t - \eta] \rightarrow \R$.  
By Theorem 2.8 in \cite{DNW} and Lemma \ref{lem:uniqueness} right below, the optimal limit path  
on $[0,t]$ (as $n \to \8$) for the coupon collector running from the proportion $\varepsilon$ to the proportion $r$ of collected coupons in time $t$ is unique and reads (we put a tilde over $g^{\varepsilon}$ below to emphasize that the interval is shifted in time and that the terminal time is $t$ and not $t - (\eta+\eps)$): 
\begin{equation}
\label{eq:g:eps}
\tilde{g}^{\varepsilon} : [0,t] \ni s \mapsto \varepsilon + \frac{1 - \varepsilon \rhops}{\rhops} \bigl( 1 - e^{-\rhops s} \bigr), 
\end{equation}
which is equal to $r$ at time $t$. By concavity of $\tilde{g}^{\eps}$,
\begin{equation*}
n \tilde{g}^{\eps}(s)  \geq n \frac{s}{t} \tilde{g}^{\eps}(t) 
+ n \big(1 - \frac{s}{t} \big) \tilde{g}^{\eps}(0) 
=  n\frac{s}{t}( r - \eps)+ n \eps  , \quad 0 \leq s \leq t. 
\end{equation*}
Coming back to $g^{\eps}$, we deduce from the inequality $t - (\eps + \eta) \leq t - t \eps/r$ that
\begin{equation}
\label{eq:ld108}
\begin{split}
n {g}^{\eps}(s)  &\geq n\frac{s - \eps}{t- (\eps+\eta)}( r - \eps)+ n \eps  
\\
&\geq  n \frac{(s- \eps)r}{t} + n \eps, \quad \eps \leq s \leq t-\eta,
\end{split}
\end{equation}
which proves that $g^{\eps}$ is strictly above the linear constraint in $A''$. Optimality of 
$g^{\eps}$ then says that (compare with 
\eqref{eq:ld20})
\begin{equation}
\label{eq:seg_opt}
\limsup_{n \rightarrow \infty } n^{-1} \ln \
\IP\big(  \hat{A},
\exists s \in [\eps,t-\eta] : N_{n}(ns) \leq n (g^{\eps}(s) - \eta')
\vert A' \bigr) < - J^\eps_{t - (\eta+\eps)}(r).
\end{equation}
Following the proof of \eqref{eq:ld105}, we deduce from \eqref{eq:ld108} and \eqref{eq:seg_opt} that \eqref{eq:ld106} holds. From (\ref{eq:ld12},\ref{eq:ld:101},\ref{eq:ld:102},\ref{eq:ld106}), we deduce that 
\begin{equation}
\label{eq:ld107}
\IP\big( A \cap A' \cap A'' \bigr) 
\geq  \exp \{-n [J_{t}(r) + \delta(\eta) + \delta'(\eps;\eta)]
+o(n) \},
\end{equation}
with $\delta(\eta) \rightarrow 0$ as $\eta \to 0$ and
$\delta'(\eps;\eta) \to 0$ as $\eps \to 0$ for a given $\eta >0$.  
\vspace{5pt}

\textit{Third Step.} We now provide a lower bound for $\IP(B \cap B' \cap B'')$ in 
\eqref{eq:ld11}. To this end, we shall use the stopping time $T_{m,n} = \inf \{ t \geq 0 : N_{n}(t) = m\}$
together with the set $C = \{\tau_{n} \geq T_{m,n} \}$. On $C$, it holds $N_{n}(\tau_{n}) \geq m$, that is $m \leq  \card \cTo(\infty)$ 
(see \eqref{def:ki} for the notations). 
Therefore, on $C$, $R(m)$ has the form:
\begin{equation*}
R(m) = \sum_{i=1}^m K(X(T_{i,n})).
\end{equation*}
In particular, on $C$,  
\begin{equation*}
R(m) \geq n \varepsilon \Rightarrow Z_{\rm tot}^{\rm GW} \geq n \varepsilon,
\end{equation*}
which implies that $
C \cap B' = C \cap 
\{R(m) \in [n \varepsilon,2n \varepsilon]\}$. Since $C$ may be also expressed as
$C = \{R(N_{n}(\ell)) > \ell; \ell =0, \ldots,  T_{m,n}-1\}$ and $\ell < T_{m,n} \Rightarrow N_{n}(\ell) \leq 
m$, we deduce that
\begin{equation}
\label{eq:measurability}
C \cap B'
\in \sigma(T_{1,n},\dots,T_{n,n}, K_{1},\dots, K_{m}).
\end{equation}
By Proposition \ref{prop:eqmod}, we deduce that $C \cap B'$ is independent of $B \cap B^{\prime \prime}$, so that 
\begin{equation}
\label{eq:ld15}
\IP \big(  B \cap B' \cap B'' \big) \geq  \IP \big(  C \cap B' \big) \IP \big(  B  \cap B'' \big).
\end{equation}
Now, we emphasize that $\tau_{n} = N_{n}(\tau_{n})$ on $C^{\complement}$ (which means that the capital is exhausted at some time less than or equal to $n$). Thus, on $C^{\complement}$,
\begin{equation*}
\tau_{n} = N_{n}(\tau_{n}) \leq N_{n}(T_{m,n}) = m \leq n \varepsilon. 
\end{equation*}
Therefore, by Lemma \ref{lem:tau}, for $\varepsilon < \varepsilon_{0}$, 
\begin{equation*}
\lim_{n \rightarrow \infty} \IP \bigl( C^{\complement} \cap \SGW \bigr) = 0. 
\end{equation*}
Finally, since $\IP(E \cap F)\geq \IP(E) - \IP(F^c)$, 
\begin{equation}
\label{eq:end:thirdstep}
\begin{split}
\IP \bigl( C \cap B' \bigr) &\geq \IP \bigl( C \cap \SGW \cap \{R(m) \in [n \varepsilon, 2n \varepsilon
] \} \bigr)
\\
&\geq \IP( C \cap \SGW ) - \IP \bigl( R(m) \not \in [n \varepsilon,2 n \varepsilon
]  \bigr)
\\
&= \IP (\SGW) - \IP \bigl( C^{\complement} \cap \SGW \bigr) - \IP \bigl( R(m) \not \in [n \varepsilon,2 n \varepsilon
]  \bigr),
\end{split}
\end{equation}
the second term in the last line converging to $0$ as $n$ tends to $\infty$. 
\vspace{5pt}

\textit{Fourth Step.} We now complete the proof when $\E K >1$, which is a simpler case to handle than the opposite case
$\E K \leq 1$.
We then choose 
\begin{equation*}
m = \lfloor \frac{\zeta}{\E K} n \varepsilon \rfloor, \quad {\rm with} \ \zeta
= \min \bigl( \frac{1+ \E K}{2}, 3/2 \bigr). 
\end{equation*} 
Clearly, $m \in [ n \varepsilon / \E K,  (1+ \E K)/ (2 \E K) \varepsilon] \subset [n \varepsilon',n (\varepsilon- \varepsilon')]$ as required, for some well-chosen $\varepsilon'$. Moreover, the typical values of 
$R(m)$ are in the neighborhood of $\zeta n \varepsilon \in (n \varepsilon,2 n \varepsilon)$. Therefore, 
by the law of large numbers, the third term in the right-hand side in \eqref{eq:end:thirdstep} tends to $0$ as $n \rightarrow \infty$. We deduce
\begin{equation}
\label{eq:ld18}
\liminf_{n \rightarrow \infty} \IP \bigl( C \cap B' \bigr) \geq \IP( \SGW ) > 0,
\end{equation}
since $\E K >1$. 
Moreover, by Mogulskii's Theorem (see Theorem 5.1.2 in \cite{dembo:zeitouni}), we have a lower bound 
for the probability that the process $(n^{-1}R(n s))_{0 \leq s \leq r+\varepsilon}$ is in the neighborhood of the path
$[0,r+\varepsilon] \ni s \mapsto s t/r$. We then observe that 
\begin{equation*}
B \cap B'' \supset
\bigl\{ (\ell-m)t/r - n\eta' \leq R(\ell) - R(m) \leq (\ell-m)t/r +  n \eta/4; m \leq \ell \leq n r + n \eps \bigr\},
\end{equation*}
since the upper condition in $B$ can be reformulated as
\begin{equation*}
R (n r + n \eps) - R(m) \leq (nr + n\eps - m) t/r + n (\eta /2- t \eps/r) + m t /r, 
\end{equation*}
and $\eta/2 - t \eps/r \geq \eta/4$. By Mogulskii's Theorem, we get:
\begin{equation}
\label{eq:ld19}
\IP( B \cap B^{\prime \prime}) \geq \exp \{-n( r+ \eps)  I(t/r)+o(n)\}.
\end{equation}
Collecting (\ref{eq:ld11},\ref{eq:ld107},\ref{eq:ld15},\ref{eq:ld18},\ref{eq:ld19}), we deduce that (for a possibly new choice of 
$\delta'(\varepsilon;\eta)$)
\begin{equation*}
\begin{split}
&\lim_{n \to \8} n^{-1} \ln
\IP\bigl(   N_n(\tau_n)/n \in [r-\eta, r+\eta], \tau_n/n \in [t-\eta, t+2\eta] \bigr) 
\\
&\hspace{15pt}\geq 
- \bigl[ J_{t}(r) + r I(t/r) \bigr] - \delta(\eta) - \delta'(\eps;\eta),
\end{split}
\end{equation*}
which tends to $-[ J_{t}(r) + r I(t/r)] - \delta(\eta)$ as $\eps \searrow 0$. Since $\eta$ can be chosen as small as needed, we complete the proof of 
\eqref{eq:lower:bound:0} in the case $\E K > 1$. 
\vspace{5pt}

\textit{Fifth Step.}
We now investigate the case $\E K \leq 1$. The above argument fails since 
$m$ cannot be chosen as $\lfloor \zeta n \varepsilon / \E K \rfloor$ on the one hand and
since $\IP(\SGW)=0$ on the other hand. 
We are thus to give a relevant version of the previous step. 
We tackle first the case $\IP(K \geq 2) >0$. The point is to change the probability measure in order to switch back to the case $\E K > 1$. The change of probability relies on the same trick as in the proof of Cramer's theorem. Since $k_*=0$ and $k^* \geq 2$, it is standard that,
for any $\xi \in (1,2)$, there exists 
$\alpha \in \R$ such that $\varphi'(\alpha)= \xi$, where
\begin{equation*}
\varphi(\alpha) = \ln \bigl[ {\mathbb E} \exp(\alpha K) \bigr]. 
\end{equation*}
We now choose $\xi = 1+ \varepsilon$ and set $\zeta'=(1+\xi)/2 \in (1, 3/2)$. For $m = \lfloor \zeta' n \varepsilon/\xi \rfloor$
($m$ is in $[n \eps',n(\eps-\eps')]$ for some well-chosen $\eps'$ and plays 
 below the same role as in the first step), we then define $\tilde{\IP}$ as
\begin{equation*}
\frac{d \tilde{\IP}}{d \IP} = \exp \biggl( \sum_{i=1}^m \{ \alpha K_{i} - \varphi(\alpha)
\}
\biggr). 
\end{equation*}
It is plain to see that, under $\tilde{\IP}$,  the variables $(K_{i})_{i \geq 1}$ are independent, the
variables $(K_{i})_{i \geq m+1}$ having the same distribution as they have under $\IP$ and the variables
$(K_{i})_{1 \leq i \leq m}$ being identically distributed with
\begin{equation*}
\tilde{\E} K_{i} = \xi, \quad i=1,\dots,m,  
\end{equation*}
where $\tilde{\E}$ denotes the expectation under $\tilde{\IP}$. Following \eqref{eq:ld15}, we write
\begin{equation}
\label{eq:ldp110}
\begin{split}
\IP(C \cap B') &= \tilde{\E} \bigl[ \frac{d \IP}{d \tilde{\IP}} {\mathbf 1}_{C \cap B'} \bigr]
\\
&= \tilde{\E} \bigl[ \exp \bigl( - \alpha R(m) +  m \varphi(\alpha) \bigr)
{\mathbf 1}_{C \cap B'} \bigr] 
\\
&\geq \exp \bigl( -  2 n \alpha \varepsilon + m \varphi(\alpha)\bigr) 
\tilde{\IP} ( C \cap B').    
\end{split}
\end{equation}
By \eqref{eq:measurability}, the probability of $C \cap B'$ under $\tilde{\IP}$ coincides with the probability 
of $C \cap B'$ under $\IP$ when the expectation of the reproduction law of the Galton-Watson tree is strictly larger than 1. Since $m = \lfloor \zeta' n \varepsilon / \tilde{\E} K \rfloor$, with $\zeta' = \min [(1+ \tilde{\E} K)/2,3/2]$, we know from the case $\E K >1$ that 
\begin{equation*}
\lim_{n \rightarrow \infty} \tilde{\IP} ( C \cap B') = \tilde{\sigma} >0,
\end{equation*}
where $\tilde{\sigma}$ stands for the probability that the Galton-Watson tree survives when the reproduction is governed by the law of $K$ under $\tilde{\IP}$. Finally, we get from \eqref{eq:ldp110}:
\begin{equation}
\label{eq:min:fifthstep} 
\liminf_{n \rightarrow \infty}
n^{-1} \ln \IP \bigl( C \cap B' \bigr) \geq - \varepsilon \xi^{-1} \bigl( 
2 \alpha \xi - \zeta' \varphi(\alpha) \bigr) =  \varepsilon G(\xi),
\end{equation}
where $G(\xi)$ remains bounded as $\xi$ ranges over any compact subset of $[1,2)$. In particular, the right hand side in \eqref{eq:min:fifthstep} tends to $0$ with $\varepsilon$. The end of the proof is then the same as in the case when 
$\IP(\SGW)>0$.

In the case $\IP(K \leq 1) = 1$, it holds $I(t/r)=\infty$ and thus ${\mathcal F}(t,r)=\infty$ as well, so that the lower bound is obvious. 
\vspace{5pt}

\textit{Final Step.}
Now we prove the lower bound on the boundary. For $r=t=0$,
we introduce $k_0=\min\{k \geq 0: \IP(K=k)>0\}$, and we write
$$
\IP \bigl(   N_n(\tau_n)/n \leq \eta, \tau_n/n \leq \eta \bigr) \geq \IP( K_1= k_0, N_n(k_0)=1)=\IP(K=k_0) n^{-k_0},
$$
showing the bound with $\cF(0,0)=0$. 

It remains to tackle the cases $0 < r=t \leq 1$ and $0 < r = 1 < t$.  Without any loss of generality, we can assume that 
$I(t/r) < \infty$ as otherwise the bound is obvious. 
Given an open set $O \ni (r,t)$ ($O \not \ni (0,0)$), we deduce from Lemma \ref{lem:continuity} 
that ${\mathcal F}$ is continuous on $O \cap \{(r',t') \in \R^2 : 0 < r ' \leq t' \wedge 1,   t'/r' \in \textrm{Dom}(I)\}$. 

If $0 < r=1 < t$, then  we can find a sequence $(r_{n},t_{n})_{n \geq 1}$, converging towards $(r,t)$, such that 
$r_{n} \nearrow r$, $t_{n} \nearrow t$ and $t_{n}/r_{n} = t/r$, with $r_{n} < 1 < t_{n}$ for any $n \geq 1$.  
Thus, 
\begin{equation}
\label{eq:finalstep}
{\mathcal F}(r,t) \geq \inf \bigl\{ {\mathcal F}(r',t'), \ (r',t') \in O' \bigr\},
\quad O' = O \cap \bigl\{(r',t') \in \R^2 : 0 < r ' < t' \wedge 1 \bigr\}. 
\end{equation}
Assuming without any loss of generality that $O \subset \{(r',t') \in \R^2 : 1- \eta < r' < 1+ \eta < t - \eta < t ' < t + \eta 
\}$, for some $\eta >0$, we deduce that 
\begin{equation}
\label{eq:finalstep:2}
\begin{split}
\liminf_{n \rightarrow \infty} n^{-1} \ln \IP \bigl( (N_{n}(\tau_{n})/n,\tau_{n}/n) \in O \bigr) 
&\geq \liminf_{n \rightarrow \infty} n^{-1} \ln \IP \bigl( (N_{n}(\tau_{n})/n,\tau_{n}/n) \in O' \bigr) 
\\
&\geq - \inf \bigl\{ {\mathcal F}(r',t'), \ (r',t') \in O' \bigr\}
\\
&= - \inf \bigl\{ {\mathcal F}(r',t'), \ (r',t') \in O \bigr\},
\end{split}
\end{equation}
the second line following from the LDP we proved above for $(r,t)$ satisfying $0 < r < t\wedge 1$. 
This is enough to conclude.

Assume $0 < r=t \leq 1$, and $I(1)<\8$ since otherwise the bound is obvious.  If
in addition $k^* \geq 2$,  \eqref{eq:finalstep} still holds,
and we can repeat \eqref{eq:finalstep:2}. 

It thus remains to handle the case when $K$ has a Bernoulli distribution of parameter $p \in (0,1]$ (if $p=0$, $K \equiv 0$ and the lower bound is obvious) and $O$ is an open set containing some point 
$(r,t)$ with $0 < r =t \leq 1$. As above, $O$ might intersect the line $(1,t')$ for $t'$ in the neighborhood of $1$; thanks to \eqref{eq:finalstep:2}, this has no real consequences.  Then, for some small $\eta >0$,
\begin{equation}
\label{eq:finalstep:3}
\begin{split}
&\IP \bigl( N_{n}(\tau_{n}) \in [n (t\!-\!\eta),n(t\!+\!\eta)], \tau_{n} \in [ n ( t\!-\!\eta, t \!+\! \eta)] \bigr)
\\
&\geq \IP \bigl( N_{n}( \lfloor n (t - \eta) \rfloor) = \lfloor n (t - \eta) \rfloor  + 1, 
 N_{n}( \lfloor n (t + \eta) \rfloor ) \leq \lfloor n (t + \eta) \rfloor
 \bigr)
\\
&\hspace{15pt} \times
\IP \bigl( R( \lfloor n(t-\eta) \rfloor +1) = \lfloor n (t-\eta) \rfloor +1   \bigr).
\end{split}
\end{equation}
Clearly,
\begin{equation}
\label{eq:finalstep:4}
\IP \bigl(   R( \lfloor n(t\!-\!\eta) \rfloor \!+\!1) = \lfloor n (t\!-\!\eta) \rfloor \!+\!1 \bigr)
= p^{\lfloor n(t\!-\!\eta) \rfloor \!+\!1}\geq p^{n(t-\eta)} 
= \exp \bigl[- n (t\!-\!\eta) I (1)  \bigr],
\end{equation}
as $I(x) = x \ln(x/p) + (1\!-\!x) \ln [(1\!-\!x)/(1\!-\!p)]$, for $x \in [0,1]$, in the Bernoulli case.  
Moreover,  
\begin{equation}
\label{eq:finalstep:5}
\begin{split}
&\IP \bigl( N_{n}( \lfloor n (t - \eta) \rfloor ) = \lfloor n (t - \eta) \rfloor  + 1, 
 N_{n}( \lfloor n (t + \eta) \rfloor) \leq \lfloor n (t + \eta) \rfloor
 \bigr)
 \\
 &\geq   ( t- \eta - 1/n ) \prod_{i=1}^{\lfloor n (t - \eta) \rfloor} (1- i/n)
 = ( t- \eta -1/n) \exp \biggl( \sum_{i=1}^{\lfloor n (t - \eta) \rfloor} \ln (1 - i/n) \biggr).
\end{split}
\end{equation}
It remains to see that 
\begin{equation}
\label{eq:finalstep:6}
\begin{split}
\lim_{n \rightarrow \infty}
- n^{-1} \sum_{i=1}^{\lfloor n (t - \eta) \rfloor} \ln (1 - i/n) 
&= - \int_{0}^{t-\eta} \ln(1-u) du 
=  (1+ \eta -t ) \ln(1+ \eta -t ) + t - \eta
\\
&\longrightarrow (1-t) \ln (1-t) + t = J_{t}(t) \qquad \textrm{as} \ \eta \rightarrow 0,
\end{split}
\end{equation}
using $\rho(t,t-\delta) \sim 2 t^{-2} \delta$ for $t \leq 1$ and $\delta \searrow 0$ in \eqref{eq:Jt(r)}. By (\ref{eq:finalstep:3},\ref{eq:finalstep:4},\ref{eq:finalstep:5},\ref{eq:finalstep:6}), the proof is easily completed. 
\qed
\vspace{5pt}

In the second step of the previous proof, we used the following.
\begin{lem}
\label{lem:uniqueness}
The path $g$ in \eqref{eq:02} is the only optimal path minimizing the limit cost for getting $r$ as proportion of collected coupons  over the rescaled time interval $[0,t]$. 
\end{lem}

\noindent  {\bf Proof:} The proof is an adaptation of Subsection A.4 in \cite{DNW}. It is sufficient to prove that the function $\varepsilon \mapsto G[\varepsilon]$ therein (see also the expression right below) is strictly convex in the neighborhood of $0_{+}$ whenever $\tilde{\gamma}$ is different from $\gamma$, where $\gamma = 1-g$ and $\tilde{\gamma}$ stands for the proportion of non-collected coupons along another path with the same boundary conditions as $\gamma$ at times $0$ and $t$.  We notice that $G[\varepsilon]$ has the form (see Subsection A.1.2 in \cite{DNW}):
\begin{equation*}
G[\varepsilon] = \int_{0}^t \bigl[ - \bigl( \dot{\gamma}_{s} + \varepsilon \dot{\eta}_{s} \bigr)
\ln \bigl( - \frac{\dot{\gamma}_{s} + \varepsilon \dot{\eta}_{s}}{\gamma_{s} + \varepsilon \eta_{s}}
\bigr) 
+ \bigl( 1+ \dot{\gamma}_{s} + \varepsilon \dot{\eta}_{s} \bigr) \ln \bigl( 
\frac{1 + \dot{\gamma}_{s} + \varepsilon \dot{\eta}_{s}}{1 - (\gamma_{s} + \varepsilon \eta_{s})} \bigr) \bigr] ds,
\end{equation*}
where $\eta = \tilde{\gamma} - \gamma$. Therefore, $G[\varepsilon]$ can be splitted into three terms:
\begin{equation*}
\begin{split}
G[\varepsilon] &= \int_{0}^t \bigl[ - \bigl( \dot{\gamma}_{s} + \varepsilon \dot{\eta}_{s} \bigr)
\ln \bigl( - \dot{\gamma}_{s} - \varepsilon \dot{\eta}_{s}
\bigr) 
+ \bigl( 1+ \dot{\gamma}_{s} + \varepsilon \dot{\eta}_{s} \bigr) \ln \bigl(
1 + \dot{\gamma}_{s} + \varepsilon \dot{\eta}_{s} \bigr) \bigr] ds
\\
&\hspace{15pt}
+ \int_{0}^t \bigl[  \bigl( \dot{\gamma}_{s} + \varepsilon \dot{\eta}_{s} \bigr)
\ln \bigl( {\gamma}_{s} + \varepsilon \eta_{s}
\bigr) 
- \bigl( \dot{\gamma}_{s} + \varepsilon \dot{\eta}_{s} \bigr) \ln \bigl( 1 - (
\gamma_{s} + \varepsilon  \eta_{s}) \bigr) \bigr] ds
\\
&\hspace{15pt}
- \int_{0}^t  \ln \bigl( 1 - (
\gamma_{s} + \varepsilon  \eta_{s}) \bigr) ds
\\
&= G_{1}[\varepsilon] + G_{2}[\varepsilon] + G_{3}[\varepsilon]. 
\end{split}
\end{equation*}
 It is well seen that $G_{1}''[\varepsilon]$ is well-defined in the neighborhood of $0$ and is always (strictly) positive unless 
 $\dot{\eta} \equiv 0$. Similarly, $G_{3}''[\varepsilon]$ is non-negative in the neighborhood of $0$. Finally, 
 \begin{equation*}
 G_{2}[\varepsilon] = \biggl[  \bigl( \gamma_{s} + \varepsilon \eta_{s} \bigr)
\ln \bigl( {\gamma}_{s} + \varepsilon \eta_{s}
\bigr) + \bigl[ 1 -   \bigl( \gamma_{s} + \varepsilon \eta_{s} \bigr) \bigr]
 \ln \bigl[ 1 - (
\gamma_{s} + \varepsilon  \eta_{s})  \bigr] \biggr]_{0}^t 
 \end{equation*}  
 is independent of $\varepsilon$ as $\eta_{t}=\eta_{0}=0$. 
\qed
\vspace{5pt}

Before proving Theorem \ref{thm:decay}, we deduce the LDP for the sequence $(N_{n}(\tau_{n}/n))_{n \geq 1}$  as an application of Theorem \ref{th:ldp}.

\begin{prop}
With the same notations as in the statement of Theorem \ref{thm:decay},  the sequence  $(N_n(\tau_n)/n)_{n \geq 1}$ obeys a LDP with rate function:
\begin{equation*}
\begin{split}
{\mathcal G}(r)  &= 
(1-r) \ln (1-r) 
\\
&\hspace{15pt}+ r \inf_{s \geq 0} \bigl\{ I \bigl(\lambda(s) \bigr) + (\lambda(s)-1) \ln \bigl( (1- \exp(-s))/r \bigr) + \lambda(s) \exp(-s)\bigr\}, 
\end{split}
\end{equation*}
if $r \in [0,1]$, and ${\mathcal G}(r)= \infty$ otherwise, 
and speed $n$. 
\end{prop}

\noindent  {\bf Proof:} By Varadhan contraction principle (see \cite[Theorem 4.2.1]{dembo:zeitouni}), we know that 
\begin{equation*}
{\mathcal G}(r) = \inf \bigl\{ {\mathcal F}(r,t), t \geq 0 \bigr\}.
\end{equation*}
Therefore,  ${\mathcal G}(r)=0$ if $r=0$, and ${\mathcal G}(r)=\infty$ if $r>1$,  as announced.
In the case $r \in (0,1]$, the infimum above can be restricted to the values of $t$ in $[r,\infty)$. Then, ${\mathcal F}(r,t)$ reads
\begin{equation}
\label{eq:contraction:1}
{\mathcal F}(r,t) = (1-r) \ln (1-r) + r \bigl[ I(\lambda) + (\lambda-1) \ln (t \rho  / \lambda r)   + \lambda \exp( - t \rho) \bigr],
\end{equation}
with $\lambda = t/r$, $t \rho$ solving the equation
$
(1 - \exp( - t \rho) )/(t \rho) = \lambda^{-1}.$ 
If $\lambda=1$, then $(\lambda-1) \ln (t \rho / \lambda r)$ and $\rho$ 
in \eqref{eq:contraction:1} are both considered as $0$. 
Letting $s= t \rho$, we note that $s$ is the unique root of the equation:
\begin{equation*}
\frac{1 - \exp( - s) }{s} = \lambda^{-1}, 
\end{equation*}
with $\lambda=1$ if $s=0$.
As $\lambda$ ranges over $[1,\infty]$, $R$ ranges over $[0,\infty)$. Expressing $\lambda$ in terms of $s$, 
the proof is easily completed.
\qed
\vspace{5pt}

\noindent
{\bf Proof  of Theorem \ref{thm:decay}:}
The LDP for $N_{n}(\tau_{n})/n$ yields
\begin{equation*}
\limsup_{n \rightarrow +  \infty} n^{-1} \ln \IP \bigl( N_{n}(\tau_{n}) = n \bigr)  \leq - {\mathcal G}(1),
\end{equation*}
with 
\begin{equation*}
{\mathcal G}(1) =
\inf_{s \geq 0} \bigl\{ I \bigl(\lambda(s) \bigr) + (\lambda(s)-1) \ln \bigl( 1 - \exp(-s)  \bigr) + \lambda(s) \exp(-s)\bigr\}.
\end{equation*}
The point is thus to prove the lower bound, which cannot be proved from the LDP directly since the lower bound in the LDP holds for open subsets only. Then, we can focus on the case $\IP(K \geq 1) >0$ as otherwise both sides in the statement of Theorem \ref{thm:decay} are infinite. We then follow the proof of the lower bound in the proof of Theorem \ref{th:ldp}. With the same notation as in the first step of the proof (in particular, given $0 < r < 1 \wedge t$), we already know that 
\begin{equation*}
\liminf_{n \rightarrow \infty}
n^{-1} \ln \IP \bigl( A \cap A' \cap A'' \cap B \cap B' \cap B'' \bigr) 
\geq - {\mathcal F}(r,t) + \delta(\eta) + \delta'(\varepsilon;\eta),
\end{equation*}
where $\delta(\eta) \rightarrow 0$ as $\eta$ tends to $0$ and
$\delta'(\varepsilon;\eta) \rightarrow 0$ as $\varepsilon$ tends to $0$ when $\eta >0$ is given. (As $r$ and $t$ do, $\varepsilon$ and $\eta$ play the same role as in the first step of the proof of the lower bound in Theorem \ref{th:ldp}.)
Define now the new events:
\begin{equation*}
\begin{split}
&A''' = \{ N_{n}(\ell +1 ) = N_{n}(\ell) + 1; \lfloor n(t-\eta) \rfloor \leq \ell <  T_{n,n} \}, 
\\ 
&B''' = \{ K_{\ell} \geq 1; \lfloor n (r+\eps) \rfloor + 1 \leq \ell \leq n \}.
\end{split}
\end{equation*} 
We claim that, on $A' \cap A'' \cap A''' \cap B' \cap B'' \cap B'''$, it holds $T_{n,n} \leq \tau_{n}$. Indeed, by definition of $A''$, 
$\lfloor n(t-\eta) \rfloor \leq T_{\lfloor n r + n \eps \rfloor,n}$, so that, 
for $\ell \in [\lfloor n (t - \eta) \rfloor+1,T_{\lfloor n r + n \eps \rfloor+1,n}-1]$, we have (in the same way as in \eqref{eq:ld:200})
\begin{equation*}
\begin{split}
R(N_{n}(\ell)) &\geq R(m) + \bigl(N_{n}(\ell)-m\bigr)t/r - n \eta'
\\
&\geq R(m) + \bigl(N_{n}(\lfloor nt - n\eta \rfloor) - m\bigr) t/r - n \eta' + (\ell - \lfloor nt - n\eta \rfloor)t/r
\\
&> \lfloor nt - n \eta \rfloor + (\ell - \lfloor nt - n \eta \rfloor) t /r > \ell,
\end{split}
\end{equation*}
and thus, for $\ell \in [T_{\lfloor n r + n \eps \rfloor+1,n},T_{n,n}]$, we also have
\begin{equation*}
\begin{split}
R(N_{n}(\ell)) &= R(\lfloor n r + n \varepsilon \rfloor) + R(N_{n}(\ell)) - R(\lfloor n r + n \varepsilon\rfloor)
\\
&= R(N_{n}(T_{\lfloor n r + n \varepsilon \rfloor,n})) + R(N_{n}(\ell)) - R(N_{n}(T_{\lfloor n r + n \varepsilon \rfloor,n}))
\\
&\geq R(N_{n}(T_{\lfloor n r + n \varepsilon \rfloor,n})) + N_{n}(\ell) -
N_{n}(T_{\lfloor n r + n \varepsilon \rfloor,n}) 
\\
&>  T_{\lfloor n r + n \eps \rfloor,n} + \ell -  T_{\lfloor n r + n \eps \rfloor,n} = \ell.
\end{split}
\end{equation*}
Thus, 
\begin{equation}
\label{eq:decay:2}
\begin{split}
\IP ( N_{n}(\tau_{n})=n ) &\geq \IP \bigl( A' \cap A'' \cap A''' \cap B' \cap B'' \cap B''' \bigr)
\\
&\geq \IP \bigl( A' \cap A'' \cap A''' \bigr) \IP \bigl( B' \cap B'' \cap B''' \bigr)
\\
&\geq \IP \bigl( A' \cap A'' \cap A''' \bigr) \IP \bigl( C \cap B' \cap B'' \cap B''' \bigr), 
\end{split}
\end{equation}
with $C$ as in the third step of the proof of the lower bound 
in Theorem \ref{th:ldp}.

From \eqref{eq:ldp:300}, it is plain to see that 
\begin{equation*}
\IP \bigl( A'''  \vert A' \cap A''  \bigr) 
\geq \exp \biggl[ \sum_{i= \lfloor n( r - \eta) \rfloor}^{\infty} \ln ( 1 - i/n) \biggr],
\end{equation*}
so that
\begin{equation}
\label{eq:decay:3}
\liminf_{n \rightarrow \infty} n^{-1} \ln \IP \bigl( A'''  \vert A' \cap A''  \bigr) 
\geq \int_{r-\eta}^1 \ln(1-u) du
= \delta''(r,\eta),
\end{equation}
where here and below $\delta''(r,\eta)$ stands for a generic term such that $\delta''(r,\eta) \rightarrow 0$ as $(r,\eta) \rightarrow (1,0)$. 
Similarly, by \eqref{eq:measurability}, 
\begin{equation*}
\IP \bigl( B''' \vert B' \cap B'' \cap C \bigr)
\geq \IP\{K \geq 1\}^{n - \lfloor n (r + \eps) \rfloor},
\end{equation*}
so that 
\begin{equation}
\label{eq:decay:4}
\liminf_{n \rightarrow \infty} n^{-1} \ln \IP \bigl(  B''' \vert C' \cap B' \cap B'' \bigr) 
\geq  (1-r) \ln  \IP(K \geq 1) = \delta'''(r),
\end{equation}
with $\delta'''(r) \rightarrow 0$ as $r \rightarrow 1$. Therefore, from (\ref{eq:decay:2},\ref{eq:decay:3},\ref{eq:decay:4}),
\begin{equation*}
\begin{split}
\liminf_{n \rightarrow \infty} n^{-1} \ln  \IP \bigl( N_{n}(\tau_{n})=n \bigr)
&\geq \liminf_{n \rightarrow \infty} n^{-1} \bigl[ \ln \IP ( A' \cap A'' ) + \ln \IP(C \cap B' \cap B'') \bigr]
\\
&\hspace{15pt}
+ \delta''(r,\eta) + \delta'''(r).
\end{split}
\end{equation*}
By the proof of the lower bound in Theorem \ref{th:ldp}, we know that 
\begin{equation*}
\liminf_{n \rightarrow \infty} n^{-1} \bigl[ \ln \IP ( A' \cap A'' ) + \ln \IP(C \cap B' \cap B'') \bigr]
\geq - {\mathcal F}(r,t) + \delta(\eta) + \delta'(\varepsilon;\eta).
\end{equation*}
In the end we deduce that 
\begin{equation*}
\liminf_{n \rightarrow \infty} n^{-1} \ln \IP \bigl( N_{n}(\tau_{n})=n \bigr)
\geq - {\mathcal F}(r,t) 
+ \delta(\eta) + \delta'(\varepsilon;\eta) + \delta''(r,\eta),
\end{equation*}
for $0 < r < t \wedge 1$. Assume then that $t >1 $ and $t \in  \textrm{Dom}^{\circ}(I)$ (interior of $\textrm{Dom}(I)$). Then, by continuity of ${\mathcal F}$ (see Lemma \ref{lem:continuity}), we can let $r$ tend to $1$.
Letting $\eta$ and $\eps$ also tend to 0, we get:
\begin{equation*}
\liminf_{n \rightarrow \infty} n^{-1} \ln \IP \bigl( N_{n}(\tau_{n})=n \bigr)
\geq - {\mathcal F}(1,t). 
\end{equation*}
If $t \geq 1$ and $t \in \textrm{Dom}(I)$, the above inequality still holds, by continuity as well, provided 
$\textrm{Dom}^{\circ}(I) \cap (1,\infty) \not = \emptyset$. If $t \not \in \textrm{Dom}(I)$, the result obviously holds, so that, in the case when 
$\textrm{Dom}^{\circ}(I) \cap (1,\infty) \not = \emptyset$, 
\begin{equation*}
\liminf_{n \rightarrow \infty} n^{-1} \ln \IP \bigl( N_{n}(\tau_{n})=n \bigr)
\geq - \inf_{t \geq 1} {\mathcal F}(1,t) = - {\mathcal G}(1).  
\end{equation*}

It thus remains to tackle the case when $\textrm{Dom}^{\circ}(I) \cap (1,\infty) = \emptyset$. Actually, this is the case when $K$ is a Bernoulli random variable. Then, we can follow the special case we discussed in the final step of the proof of Theorem \ref{th:ldp}. Indeed, for $0< r < 1$, we deduce from  (\ref{eq:finalstep:3},\ref{eq:finalstep:4},\ref{eq:finalstep:5}) (with $t=1$ and $r=t-\eta$ therein), 
\begin{equation*}
\IP \bigl( N_{n}(\tau_{n})/n \geq r \bigr) \geq (r-1/n) \exp \biggl( - \sum_{i=1}^{\lfloor nr \rfloor}
\ln (1 - i/n) - n r I(1) \biggr). 
\end{equation*} 
Letting $r$ tend to $1$, we deduce that 
\begin{equation*}
\IP \bigl( N_{n}(\tau_{n}) = n \bigr) \geq  \exp \biggl( - \sum_{i=1}^{n-1}
\ln (1 - i/n) - n I(1) \biggr),
\end{equation*} 
so that 
\begin{equation*}
\liminf_{n \rightarrow \infty} n^{-1} \ln  
\IP \bigl( N_{n}(\tau_{n}) = n \bigr)
 \geq - \int_{0}^1 \ln(1-u) du - I(1) = - (1 + I(1)). 
\end{equation*}
It then remains to check that it is equal to $-{\mathcal G}(1)$. Clearly, the infimum in the definition of ${\mathcal G}(1)$ is reduced to the $s$'s such that $\lambda(s) =1$, that is to $s=0$. We easily deduce that ${\mathcal G}(1) = 1 + I(1)$. 

This completes the proof of the variational formula. Since it has compact level sets, 
the lower semi-continuous function $\cF$, when restricted to the set of points $(r,t)$ with $r=1$,
achieves its minimum, and the value of the minimum is non zero. Hence ${\mathcal G}(1) >0$. \qed

{\bf Acknowledgements:} The authors thank anonymous referees for their constructive comments 
and careful reading of the paper which allowed us to improve the presentation. R.S. is grateful to
G. Stacey Staples for his Mathematica expertise.  

\thispagestyle{empty}

\end{document}